\documentstyle[amscd,amssymb,verbatim,11pt]{amsart}

\theoremstyle{plain}
\newtheorem{Prop}{Proposition}[section]
\theoremstyle{definition}
\newtheorem{Def}[Prop]{Definition}
\newtheorem{Thm}[Prop]{Theorem}
\newtheorem{Cor}[Prop]{Corollary}

\newtheorem{Lem}[Prop]{Lemma}
\theoremstyle{remark}
\newtheorem{Rem}[Prop]{Remark}
\newtheorem{Problem}[Prop]{\bf Problem}

\numberwithin{equation}{section}
\begin{document}

\title[Isomorphisms in pro-categories]%
    {Isomorphisms in pro-categories}

\author{J.Dydak and F.R. Ruiz del Portal}
\date{November 3, 2003}
\keywords{isomorphism,monomorphism, epimorphism, pro-categories}

\subjclass{16B50, 18D35, 54C56}

\thanks{The first-named author supported in part by grant DMS-0072356
from NSF and by the Ministry of Science and
Education of Spain. The second-named author supported by MCyT, BMF2000-0804-C03-01.}

\begin{abstract}
A morphism of a category which is simultaneously an epimorphism and a 
monomorphism is
called a bimorphism. In \cite{DR2} we gave characterizations of monomorphisms
(resp. epimorphisms) in arbitrary pro-categories, pro-\(C\), where 
\(C\) has direct sums
(resp. weak push-outs). In this paper we introduce the notions of strong 
monomorphism and strong epimorphism.
Part of their significance is that they are preserved by functors.
These notions and their characterizations lead us to important 
classical properties
and problems in shape and pro-homotopy. For instance, strong 
epimorphisms allow us
to give a categorical point of view of uniform movability
and to introduce a new kind of movability, the sequential movability.
Strong monomorphisms are connected to a problem of K.Borsuk regarding
a descending chain of retracts of ANRs.
If \(f: X \rightarrow Y\) is a bimorphism in the pointed shape 
category of topological spaces, we prove that \(f\) is a weak 
isomorphism and \(f\) is an isomorphism provided \(Y\) is 
sequentially movable and $X$ or $Y$ is the suspension of a 
topological space.
If \(f: X \rightarrow Y \) is a bimorphism
in the pro-category pro-\(H_0\) (consisting of inverse systems in 
\(H_0\), the homotopy category of pointed connected CW complexes) we 
show that \(f\) is an isomorphism provided \(Y\) is sequentially 
movable.
\end{abstract}

\maketitle

\medskip
\medskip
\tableofcontents
\section{Introduction}

\medskip

The fundamental problem in any category $C$ is to detect
its isomorphisms. A way to do it is by considering,
as in the category of groups,
  the notions of epimorphism and monomorphism in
abstract
categories.
\begin{Def}\label{XXX1.1}
\par A morphism $f:X\to Y$ of a category $ C$ is called an
{\bf epimorphism} if the induced function
$f^*:Mor(Y,Z)\to Mor(X,Z)$ is one-to-one for each object Z of
$ C$.
\par A morphism $f:X\to Y$ of a category $ C$ is called a
{\bf monomorphism} if the induced function
$f_*:Mor(Z,X)\to Mor(Z,Y)$ is one-to-one for each object Z of
$ C$.
\end{Def}

In common terms, $f$ is
an epimorphism (respectively, monomorphism) of $ C$ if
  $g\circ f=h\circ f$ (respectively, $f\circ g=f\circ h$)
implies $g=h$ for any two morphisms $g,h:Y\to Z$ (respectively,
$g,h:Z\to X$).
\par The main drawback of the two concepts is that they are not
functorial, i.e. they are not preserved by covariant functors.
One can easily check that morphisms $f$ of a category $C$
with the property that $F(f)$ is a monomorphism (respectively,
an epimorphism) of the category $D$ for any covariant
functor $F:C\to D$ are exactly those having a left (respectively, right) inverse.
Having a left (respectively, right) inverse is, obviously,
 a functorial property.

\par A well-known and easy exercise is the following.
\begin{Prop}\label{XXX1.2} A monomorphism (respectively, epimorphism)
which has a left (respectively, right) inverse
  is an isomorphism.
\end{Prop}

The main object of our study are isomorphisms in pro-categories
(see a review of pro-categories in the next section). In the case of pro-categories
one can consider the following variant of functoriality:
Suppose $f$ is a morphism of pro-$C$. When is $(pro-F)(f)$
a monomorphism (respectively, an epimorphism) of pro-$D$
for any covariant functor $F:C\to D$?
\par It turns out (see \ref{XXapp.1})
that those are exactly strong monomorphisms
(respectively, strong epimorphisms) - the key concepts for this paper
(see \ref{XXX2.2} for a definition). Our best results characterizing isomorphisms
of pro-categories are: \ref{XXX2.10} (stating that $f$
is an isomorphism of pro-$C$ if and only if it is a strong monomorphism
and an epimorphism) and \ref{XXX2.12}
(stating that, for categories $C$ with direct sums, $f$
is an isomorphism of pro-$C$ if and only if it is a strong epimorphism
and a monomorphism).
Our best general application of strong epimorphisms is a characterization
of uniform movability in  \ref{XXX3.2} with the resulting
characterization of isomorphisms $f:X\to Y$ such that
$Y$ is uniformly movable and $\underset{\leftarrow}{\text{lim}}(f)$ 
is an isomorphism of $C$
 (see  \ref{XXX3.3}).

\begin{Def}\label{XXX1.3}
\par A morphism $f:X\to Y$ of a category $ C$ is called an
{\bf bimorphism} if it is both an epimorphism and a monomorphism
of $ C$.
\par
A category $ C$ is called {\bf balanced} if every bimorphism
of $ C$ is an isomorphism.
\end{Def}

\par The following question was posed in \cite{DR}.

\par
\begin{Problem}\label{XXX1.4} Suppose a category $C$ is balanced.
Is the pro-category pro-$C$ balanced?
\end{Problem}

This question was answered negatively in \cite{DR2}, so an amended
version of it is as follows.
\begin{Problem}\label{XXX1.4.5} Suppose $C$ is balanced category
with direct sums and weak push-outs.
Is the pro-category pro-$C$ balanced?
\end{Problem}

A natural question is to decide which common categories
are balanced. It is so
  in the case of the homotopy category $H_0$ of pointed connected
CW complexes. The question of whether $H_0$
is balanced has been open for a while with Dyer and Roitberg \cite{DyR}
resolving it in positive and Dydak \cite{D2} giving a simple proof
of it. Mukherjee
\cite{Mu} generalized the approach from \cite{DyR} to the equivariant case and
Mor\'on-Ruiz del Portal \cite{MoP}
showed that the shape category of pointed, movable, metric continua
is not balanced but every weak isomorphism is a bimorphism.
We recommend \cite{G} for a near complete list and a
thorough
review of results related to monomorphisms and epimorphisms of $H_0$.

\par In \cite{DR} the authors embarked on a program to determine
if pro-$H_0$ is balanced and that paper contains results
on bimorphisms in tow($H_0$), the category of towers
in $H_0$. Section 8 of this paper generalizes
those results to bimorphisms of pro-$H_0$.
In section 9 we investigate bimorphisms of the shape category
and in section 10 we relate the concept of strong monomorphism
to a question of K.Borsuk.
\par
The authors are grateful to 
the referee for numerous improvements of the paper.
We are indebted to M.A. Mor\'on for help in understanding
of the Borsuk's problem.

\section{Review of pro-categories}

\par Let us recall basic facts about pro-categories
(for details see \cite{MS}).
Loosely speaking, the pro-category pro-$C$ of $C$
is the universal category with inverse limits containing
$C$ as a full subcategory. Quite often one considers
pro-objects indexed by small cofiltered categories.
However, those are isomorphic to pro-objects indexed
by directed sets (see pp.14-15 of \cite{MS}), so in this paper the
 objects of pro-$C$ are inverse systems
$X=(X_\alpha,p^\beta_\alpha,A)$ in $C$ such that $A$
is a directed set.
To simplify notation we will call $A$ the {\bf index set}
of $X$ and we will denote it by $I(X)$.
Given $\alpha,\beta\in I(X)$ with $\alpha<\beta$, the {\bf bonding
morphism} $p^\beta_\alpha$ from $X_\beta$ to $X_\alpha$ will be denoted
by $p(X)^\beta_\alpha$.
\par Given an inverse system $X$ in $C$ and $P\in Ob(C)$
($P$ is an object of $C$), the set of morphisms
of pro-$C$ from $X$ to $P$ is the direct limit
of $Mor(X_\alpha,P)$, $\alpha\in I(X)$.
Thus each morphism $f$ from $X$ to $P$ has
its {\bf representative} $g:X_\alpha\to P$,
and two representatives $g:X_\alpha\to P$
and $h:X_\beta\to P$ determine the same morphism
if there is $\gamma>\alpha,\beta$ with
$g\circ p(X)^\gamma_\alpha=h\circ p(X)^\gamma_\beta$.
In particular, the morphism from $X$ to $X_\alpha$
represented by the identity $X_\alpha\to X_\alpha$
is called the {\bf projection morphism}
and denoted by $p(X)_\alpha$.
It is clear how to compose morphisms from $X$ to $P$
with morphisms from $P$ to $Q$ if $P,Q\in Ob(C)$.
\par
If $X$ and $Y$ are two inverse systems in $C$,
then any morphism $f:X\to Y$ of pro-$C$
can be identified with a family of morphisms
$\{f_\alpha:X\to Y_\alpha\}_{\alpha\in I(Y)}$
such that $p(Y)^\beta_\alpha\circ f_\beta=f_\alpha$
for all $\alpha<\beta$ in $I(Y)$.
Notice that $f_\alpha=p(Y)_\alpha\circ f$ for each $\alpha\in I(Y)$.
Therefore one has a simple characterization of isomorphisms
of pro-$C$.

\begin{Prop}\label{XXX1.5} A morphism $f:X\to Y$ of pro-$C$
is an isomorphism if and only if $f^\ast:Mor(Y,P)\to Mor(X,P)$
is a bijection for all $P\in Ob(C)$.
\end{Prop}

Of particular interest are isomorphisms $f:X\to P\in Ob(C)$.
If such an isomorphism exists, then $X$ is called {\bf stable}.

\par If $s$ is a directed subset of $I(X)$,
then by $X_s$ we will denote the induced
inverse system $(X_\alpha,p(X)^\beta_\alpha,s)$.
Notice that the family $\{p_\alpha:X\to X_\alpha\}_{\alpha\in s}$
induces a morphism from $X$ to $X_s$
which will be denoted by $p(X)_s$.
If $s$ is a {\bf cofinal} subset of $I(X)$
(that means for any $\alpha\in I(X)$ there is $\beta\in s$
so that $\beta>\alpha$), then $p(X)_s$
is an isomorphism of pro-$C$.
\par Of particular use are {\bf level morphisms}
of pro-$C$. Those are $f:X\to Y$, where $X$ and $Y$ have
identical index sets and there are
representatives $f_\alpha:X_\alpha\to Y_\alpha$
of $p(Y)_\alpha\circ f$,
$\alpha\in I(X)$, such that
$p(Y)^\beta_\alpha\circ f_\beta=f_\alpha\circ p(X)^\beta_\alpha$
for all $\alpha<\beta$.
In such a case we say that $f$ is a
{\bf level morphism induced by the family}
$\{f_\alpha:X_\alpha\to Y_\alpha\}_{\alpha\in I(Y)}$.
\par It is also convenient to consider inverse systems
$X$ such that $I(X)$ is a {\bf cofinite directed set}
which means that each element of $I(X)$ has only finitely
many predecessors. The following result
is of particular use (see \cite{MS}, Theorem 3 on p.12).
\begin{Prop}\label{XXX1.6} For any morphism $f:X\to Y$ of pro-$C$
there exists a level morphism $f':X'\to Y'$
and isomorphisms $i:X\to X'$, $j:Y'\to Y$
such that $f=j\circ f'\circ i$ and $I(X')$
is a cofinite directed set.
Moreover, the bonding morphisms of $X'$ (respectively, $Y'$)
are chosen from the set of bonding morphisms of $X$
(respectively, $Y$).
\end{Prop}

In the special case of $X$ being an object of $C$
one can create $X'$ by putting $X'_\alpha=X$
and $p(X')^\beta_\alpha=id_X$ for each $\beta>\alpha$
in $I(Y)$. Notice that in this case $Y'=Y$ and
$f'$ is induced by the family $\{p(Y)_\alpha\circ f\}_{\alpha\in I(Y)}$.
In what follows morphisms from objects $X$ of $C$
to inverse systems $Y$ in $C$ will be automatically replaced
by level morphisms from $X'$ to $Y$. This is needed
as part of our strategy is to select increasing
sequences in $I(X)$ which is not possible
if $I(X)$ contains an upper bound (which implies
that $X$ is {\bf stable}, i.e. isomorphic to
an object of $C$).
\par
Another reason level morphisms are very useful is that one has
a very simple criterion of them being an isomorphism
  (see \cite{MS}, Theorem 5 on p.112).
\begin{Prop}\label{XXX1.6b} A level morphism $f:X\to Y$ of pro-$C$
is an isomorphism if and only if for each $\alpha \in I(X)$
there is $\beta >\alpha $ and $g:Y_\beta \to X_\alpha $
such that $f_\alpha \circ g=p(Y)^\beta _\alpha $
and $g\circ f_\beta =p(X)^\beta _\alpha $.
\end{Prop}

We will need the following characterization of monomorphisms
in pro-$C$ such that $C$ has direct sums (see \cite{DR2}, Proposition 2.3).


\begin{Prop}\label{XXX1.mono}
Suppose that \(f=\{f_{\alpha}:X_{\alpha} \rightarrow 
Y_{\alpha}\}_{\alpha \in I(X)}\) is a level morphism of pro-\(C\). 
Consider the following conditions:

a. \(f\) is a monomorphism.

b. For each \(\alpha \in I(X)\) there is \(\beta \in I(X)\), 
\(\beta > \alpha\), such that for any
\(u,v: P \in Ob(C) \rightarrow X_{\beta}\), \(f_{\beta} \circ u= 
f_{\beta} \circ v\) implies that
\(p(X)_{\alpha}^{\beta} \circ u = p(X)_{\alpha}^{\beta} \circ v\).

c. For each \(\alpha \in I(X)\) there is \(\beta \in I(X)\), 
\(\beta > \alpha\), such that for any
\(u,v: T \rightarrow X_{\beta}\), \(f_{\beta} \circ u= f_{\beta} 
\circ v\) implies that
\(p(X)_{\alpha}^{\beta} \circ u = p(X)_{\alpha}^{\beta} \circ v\).

Conditions b) and c) are equivalent and imply Condition a).
If \(C\) has direct sums, then all three conditions are equivalent
\end{Prop}

Also, we will give a characterization of epimorphisms
in pro-$C$ such that $C$ has weak push-outs.

\begin{Def}\label{XXXpushouts}
A commutative diagram
\begin{center}
\setlength{\unitlength}{1mm}
\begin{picture}(63,28)

\put(5.3,22){\(A\)}

\put(58.3,22){\(B\)}

\put(5.3,0){\(C\)}

\put(58.3,0){\(P\)}

\put(32,25){\(f\)}

\put(32,3){\(v\)}

\put(61,11){\(u\)}

\put(1,11){\(g\)}

\put(15,23){\vector(1,0){37}}

\put(15,1){\vector(1,0){37}}

\put(6.3,19){\vector(0,-1){13}}

\put(59.3,19){\vector(0,-1){13}}

\end{picture}
\end{center}
in category $C$ is a {\bf weak push-out} 
(respectively, {\bf push-out}) of the diagram
\begin{center}
\setlength{\unitlength}{1mm}
\begin{picture}(63,28)

\put(5.3,22){\(A\)}

\put(58.3,22){\(B\)}

\put(5.3,0){\(C\)}

\put(32,25){\(f\)}

\put(1,11){\(g\)}

\put(15,23){\vector(1,0){37}}

\put(6.3,19){\vector(0,-1){13}}
\end{picture}
\end{center}
if for any commutative diagram
\begin{center}
\setlength{\unitlength}{1mm}
\begin{picture}(63,28)

\put(5.3,22){\(A\)}

\put(58.3,22){\(B\)}

\put(5.3,0){\(C\)}

\put(58.3,0){\(P'\)}

\put(32,25){\(f\)}

\put(32,3){\(v'\)}

\put(61,11){\(u'\)}

\put(1,11){\(g\)}

\put(15,23){\vector(1,0){37}}

\put(15,1){\vector(1,0){37}}

\put(6.3,19){\vector(0,-1){13}}

\put(59.3,19){\vector(0,-1){13}}

\end{picture}
\end{center}
in $C$ there is a morphism (respectively, a unique morphism)
$t:P\to P'$ such that 
$u'=t\circ u$ and $v'=t\circ v$.
If every diagram 
\begin{center}
\setlength{\unitlength}{1mm}
\begin{picture}(63,28)

\put(5.3,22){\(A\)}

\put(58.3,22){\(B\)}

\put(5.3,0){\(C\)}

\put(32,25){\(f\)}

\put(1,11){\(g\)}

\put(15,23){\vector(1,0){37}}

\put(6.3,19){\vector(0,-1){13}}
\end{picture}
\end{center}
in $C$ has a weak push-out,
then we say that $C$ is a category with {\bf weak push-outs}.
\end{Def}


\begin{Prop}\label{XXXpushoutsInCW} The homotopy category $H_0$
of pointed connected CW complexes is a category with weak push-outs
but not a category with push-outs.
\end{Prop}

\begin{pf}
Suppose 
\begin{center}
\setlength{\unitlength}{1mm}
\begin{picture}(63,28)

\put(5.3,22){\(A\)}

\put(58.3,22){\(B\)}

\put(5.3,0){\(C\)}

\put(32,25){\(f\)}

\put(1,11){\(g\)}

\put(15,23){\vector(1,0){37}}

\put(6.3,19){\vector(0,-1){13}}
\end{picture}
\end{center}
is a commutative diagram in $H_0$.
Pick cellular maps $b:A\to B$ and $c:A\to C$
representing $f$ and $g$, respectively.
Let $P$ be the union of the reduced mapping cylinders $M(a)$ and $M(b)$
so that $M(a)\cap M(b)=A$.
There are natural inclusions $u:C\to P$ and $v:B\to P$.
Given maps $u':C\to P'$ and $v':B\to P'$ so that
$u'\circ c$ is homotopic to $v'\circ b$, any homotopy
$H$ from $u'\circ c$ to $v'\circ b$ leads naturally to a map 
$t:P\to P'$ so that $t$ extends both $u'$ and $v'$.
\par To show that $H_0$ does not have push-outs, let us use
an example provided to us by the referee. Namely, both $A$ and $C$
are the unit circle $S^1$, $B$ is trivial, and $g(z)=z^2$ for $z\in S^1$.
Suppose
\begin{center}
\setlength{\unitlength}{1mm}
\begin{picture}(63,28)

\put(5.3,22){\(A\)}

\put(58.3,22){\(B\)}

\put(5.3,0){\(C\)}

\put(58.3,0){\(P\)}

\put(32,25){\(f\)}

\put(32,3){\(v\)}

\put(61,11){\(u\)}

\put(1,11){\(g\)}

\put(15,23){\vector(1,0){37}}

\put(15,1){\vector(1,0){37}}

\put(6.3,19){\vector(0,-1){13}}

\put(59.3,19){\vector(0,-1){13}}

\end{picture}
\end{center}
is a push-out diagram.
Given any pointed CW complex $X$ and an element $\alpha\in \pi_1(X)$
satisfying $\alpha^2=1$, the diagram
\begin{center}
\setlength{\unitlength}{1mm}
\begin{picture}(63,28)

\put(5.3,22){\(A\)}

\put(58.3,22){\(B\)}

\put(5.3,0){\(C\)}

\put(58.3,0){\(X\)}

\put(32,25){\(f\)}

\put(32,3){\(\alpha\)}

\put(61,11){\(const\)}

\put(1,11){\(g\)}

\put(15,23){\vector(1,0){37}}

\put(15,1){\vector(1,0){37}}

\put(6.3,19){\vector(0,-1){13}}

\put(59.3,19){\vector(0,-1){13}}

\end{picture}
\end{center}
is commutative, so there is $h_\alpha:P\to X$
with $\alpha=h_\alpha\circ u$.
That correspondence between $\alpha$ and $h_\alpha$ establishes
a bijection between morphisms from $P$ to $X$
and the subset of $\pi_1(X)$ consisting of elements whose square is 1.
In particular, all higher cohomology of $P$ vanishes
and $H^1(P;G)$ consists of elements of $G$ order two or less.
The contradiction arises by looking at the exact sequence
$0\to Z/2\to Z/4\to Z/2\to 0$ and induced exact sequence
of cohomology groups (recall that higher cohomology of $P$ vanishes)
$0\to H^1(P;Z/2)\to H^1(P;Z/4)\to H^1(P;Z/2)\to 0$.
It contradicts that $H^1(P;G)=\{g\in G | 2g=0\}$
as there is no exact sequence $0\to Z/2\to Z/2\to Z/2\to 0$.
\end{pf}

\begin{Prop}\label{XXX1.epi} Suppose 
 \(f=\{f_{\alpha}:X_{\alpha} \rightarrow Y_{\alpha}\}_{\alpha \in I(X)}\) is a level morphism of pro-$C$ and consider the following conditions:

a. for each \(\alpha \in I(X)\) there is \(\beta \in I(X)\), \(\beta > \alpha\), such that for any
\(u,v: Y_{\alpha} \rightarrow P\in Ob(C)\), \(u \circ f_{\alpha}= v \circ f_{\alpha}\) implies that \( u \circ p(Y)_{\alpha}^{\beta}  = v \circ p(Y)_{\alpha}^{\beta}\).

b. \(f\) is an epimorphism  of pro-$C$.

Condition a) is stronger than Condition b). If \(C\) is a category with weak push-outs,
then both conditions are equivalent.

\end{Prop}

\begin{pf}
\par a)$\implies$b). It suffices to show that $u,v:Y\to P\in Ob(C)$
and $u\circ f=v\circ f$ implies $u=v$.
Pick representatives $u':Y_\alpha\to P$
of $u$ and $v':Y_\alpha\to P$
of $v$ for some $\alpha\in I(X)$.
We may assume that $u'\circ f_\alpha=v'\circ f_\alpha$.
There is $\beta>\alpha$ such that 
\( u' \circ p(Y)_{\alpha}^{\beta}  = v' \circ p(Y)_{\alpha}^{\beta}\)
which implies $u=v$.
\par b)$\implies$a) if \(C\) is a category with weak push-outs.
Let

\begin{center}
\setlength{\unitlength}{1mm}
\begin{picture}(63,28)

\put(5.3,22){\(X_{\alpha}\)}

\put(58.3,22){\(Y_{\alpha}\)}

\put(5.3,0){\(Y_{\alpha}\)}

\put(58.3,0){\(M\)}

\put(32,25){\(f_{\alpha}\)}

\put(32,3){\(b\)}

\put(61,11){\(a\)}

\put(1,11){\(f_{\alpha}\)}

\put(15,23){\vector(1,0){37}}

\put(15,1){\vector(1,0){37}}

\put(6.3,19){\vector(0,-1){13}}

\put(59.3,19){\vector(0,-1){13}}

\end{picture}
\end{center}

\noindent be a weak push-out of 

\begin{center}
\setlength{\unitlength}{1mm}
\begin{picture}(63,28)

\put(5.3,22){\(X_{\alpha}\)}

\put(58.3,22){\(Y_{\alpha}\)}

\put(5.3,0){\(Y_{\alpha}\)}

\put(32,25){\(f_{\alpha}\)}

\put(1,11){\(f_{\alpha}\)}

\put(15,23){\vector(1,0){37}}


\put(6.3,19){\vector(0,-1){13}}


\end{picture}
\end{center}
\noindent
 There is \(\beta > \alpha\) so that \(a \circ p(Y)_{\alpha}^{\beta}= b \circ p(Y)_{\alpha}^{\beta}\)
as \(f\) is an epimorphism. 
If \(u,v: Y_{\alpha} \rightarrow P\in Ob(C)\), \(u \circ f_{\alpha}= v \circ f_{\alpha}\), then

\begin{center}
\setlength{\unitlength}{1mm}
\begin{picture}(63,28)

\put(5.3,22){\(X_{\alpha}\)}

\put(58.3,22){\(Y_{\alpha}\)}

\put(5.3,0){\(Y_{\alpha}\)}

\put(58.3,0){\(P\)}

\put(32,25){\(f_{\alpha}\)}

\put(32,3){\(v\)}

\put(61,11){\(u\)}

\put(1,11){\(f_{\alpha}\)}

\put(15,23){\vector(1,0){37}}

\put(15,1){\vector(1,0){37}}

\put(6.3,19){\vector(0,-1){13}}

\put(59.3,19){\vector(0,-1){13}}

\end{picture}
\end{center}

\noindent is commutative, so there is \(i: M \rightarrow P\)  such that
\( i \circ b= v\) and \( i \circ a = u\).

Therefore,

\[u \circ p(Y)_{\alpha}^{\beta} = i \circ a \circ p(Y)_{\alpha}^{\beta} = i \circ b \circ p(Y)_{\alpha}^{\beta}= v \circ p(Y)_{\alpha}^{\beta}.\]
\end{pf}

\begin{Rem}\label{XXXrefreeThanks}
In the original version of this paper Proposition \ref{XXX1.epi}
was stated and used for categories $C$ with push-outs.
However, our main application was for $C=H_0$ and it was pointed out
by the referee that $H_0$ does not have push-outs
(see the proof of \ref{XXXpushoutsInCW}). That is how the notion of weak push-outs
was created to salvage \ref{XXX1.epi} and allow
all applications to be valid.
\end{Rem}

\par Let us proceed with an abstract construction. Its meaning will
be explained shortly.
\begin{Def}\label{XXX1.7} Suppose $c$ is an ordinal.
Given $X\in Ob(pro-C)$ an object
$Sub_c(X)$ of pro-(pro-$C$) is defined as follows:

1. the index set $I(Sub_c(X))$ of $Sub_c(X)$ consists
of all increasing functions $s:\{n\mid n < c\}\to I(X)$,
where $n$ is a cardinal number smaller than $c$,

2. $s\leq t$, $s,t\in I(Sub_c(X))$, holds if and only if
$s(n)\leq t(n)$ for all $n < c$,

3. $Sub_c(X)_s:=(X_{s(n)},p(X)^{s(m)}_{s(n)})$ for all $s\in I(Sub_c(X))$,

4. $p(Sub_c(X))_s^t$ is the level morphism
induced by $\{p(X)^{t(n)}_{s(n)}\}_{n < c}$ for all $s<t$.
\end{Def}

An important case of \ref{XXX1.7} is $c=2$:
$Sub_2(X)$ is another way of looking at $X$
as an object of pro-(pro-$C$) (the standard way
corresponds to the canonical embedding of any category into
its pro-category). Observe that the projections $Sub_c(X)_s\to X_{s(1)}$
induce a morphism from $Sub_c(X)$ to $Sub_2(X)$
for any $c\ge 2$.
Also notice that the family
$\{p(X)_\alpha:X\to X_\alpha\}_{\alpha\in I(X)}$
induces a morphism from $X$ to $Sub_2(X)$ of pro-(pro-$C$). 
That morphism will be used
later on to explain the concept of uniform movability.

\par Another important case of \ref{XXX1.7} is $c=\omega_0$
as part of our work is related to reducing properties
of pro-$C$ to the properties of its full subcategory tow($C$)
consisting of {\bf towers}, i.e. inverse sequences in $C$.

\par Just as every morphism of pro-$C$ from $X$
to an object of $C$ factors through a subterm of $X$,
every morphism from $X$ to a tower factors through
a subtower.
\begin{Prop}\label{XXX1.8} Suppose $C$ is a category.
If $f:X\to Y$ is a morphism to a tower,
then there is a subtower $X_s$ of $X$
and a level morphism $g:X_s\to Y$
such that $f=g\circ p(X)_s$.
\end{Prop}
\begin{pf} Choose $s(1)\in I(X)$ such that there is
a representative $g_1:X_{s(1)}\to Y_1$
of $p(Y)_1\circ f$.
Suppose $s(i)$ and $g_i:X_{s(i)}\to Y_i$ are defined for $i\leq n$
such that $g_i$ is a representative of $p(Y)_i\circ f$.
Find $\alpha\in I(X)$, $\alpha>s(n)$
such that there is a representative $h:X_\alpha\to Y_{n+1}$
of $p(Y)_{n+1}\circ f$.
Since both $g_n$ and $p(Y)_n^{n+1}\circ h$ are representatives
of $p(Y)_n\circ f$, there is $s(n+1)>\alpha$
such that $g_n\circ p(X)^{s(n+1)}_{s(n)}=
p(Y)_n^{n+1}\circ h\circ p(X)^{s(n+1)}_\alpha$.
By putting $g_{n+1}=h\circ p(X)^{s(n+1)}_\alpha$
we complete the inductive construction
of a level morphism $g:X_s\to Y$ satisfying
$g\circ p(X)_s=f$.
\end{pf}

\begin{Prop}\label{XXX1.9} Suppose $C$ is a category
and $f:Y\to Z$ is a morphism of towers in $C$.
If $X_s$ is a subtower of $X$ and $g,h:X_s\to Y$
are two morphisms such that $f\circ g\circ p(X)_s=
f\circ h\circ p(X)_s$,
then there is a subtower $X_t$ of $X$
such that $t>s$ and $f\circ g\circ p(X)_s^t=
f\circ h\circ p(X)_s^t$.
\end{Prop}
\begin{pf} Special Case: $f,g$, and $h$ are level morphisms.
\par For each $n\in N$ the morphisms $f_n\circ g_n$
and $f_n\circ h_n$ are representatives of
$p(Z)_n\circ f\circ g\circ p(X)_s$,
so there is $t(n)>s(n)$ such that
$f_n\circ g_n\circ p(X)^{t(n)}_{s(n)}=f_n\circ h_n\circ p(X)^{t(n)}_{s(n)}$.
Using induction one can ensure $t(n)>t(n-1)$
which completes the proof of Special Case.
\par General Case. By \ref{XXX1.8} there is an increasing sequence $u:N\to N$
and a level morphism $f':Y_u\to Z$ so that $f=f'\circ p(Y)_u$.
Using \ref{XXX1.8} again one can find an increasing sequence
$v:N\to im(s)$ and level morphisms $g',h':X_v\to Y_u$
such that $g'\circ p(X_s)_v=p(Y)_u\circ g$ and
$h'\circ p(X_s)_v=p(Y)_u\circ h$.
Since $f'\circ g'\circ p(X)_v=
f'\circ g'\circ p(X_s)_v\circ p(X)_s=
f'\circ p(Y)_u \circ g\circ p(X)_s=
f\circ g\circ p(X)_s$
and, similarly,
$f'\circ h'\circ p(X)_v=
f\circ h\circ p(X)_s$, we get
$f'\circ g'\circ p(X)_v=f'\circ h'\circ p(X)_v$.
By Special Case there is $t>v$ so that
$f'\circ g'\circ p(X)_v^t=f'\circ h'\circ p(X)_v^t$.
Now $f\circ g\circ p(X)_s^t=f'\circ p(Y)_u\circ g\circ p(X)_s^t=
f'\circ g'\circ p(X_s)_v\circ p(X)_s^t=
f'\circ g'\circ p(X)_v^t$ and, similarly,
$f\circ h\circ p(X)_s^t=
f'\circ h'\circ p(X)_v^t$. Therefore
$f\circ g\circ p(X)_s^t=f\circ h\circ p(X)_s^t$.
\end{pf}

\begin{Cor}\label{XXX1.10} Let $C$ be a category.

1. Every monomorphism (respectively, epimorphism)
of $C$ is a monomorphism (respectively, epimorphism) of pro-$C$.

2. Every monomorphism (respectively, epimorphism)
of $tow(C)$ is a monomorphism (respectively, epimorphism) of pro-$C$.

3. Every bimorphism of $C$ or $tow(C)$ is a bimorphism of pro-$C$.
\end{Cor}
\begin{pf} A. Let us prove 1) and 2) for epimorphisms.
Suppose $f:X\to Y$ is an epimorphism of $D$, $D=C$
or $D=tow(C)$, and $g,h:Y\to Z$ satisfy $g\circ f=h\circ f$.
To show $g=h$ it suffices to prove
$p(Z)_\alpha\circ g=p(Z)_\alpha\circ h$ for all
$\alpha\in I(Z)$. Since $Z_\alpha$ is an object of $D$
and $(p(Z)_\alpha\circ g)\circ f=(p(Z)_\alpha\circ h)\circ f$,
one gets $p(Z)_\alpha\circ g=p(Z)_\alpha\circ h$ as $f$
is an epimorphism of $D$.
\par B. A. Let us prove 1) and 2) for monomorphisms.
 Suppose $f:X\to Y$ is a monomorphism of $D$, $D=C$
or $D=tow(C)$, and $g,h:Z\to X`$ satisfy $f\circ g=f\circ h$.
By \ref{XXX1.8} and \ref{XXX1.9} there is a sequence $s$ in $I(X)$ 
(an element $s$ of $I(X)$
if $D=C$) and there are level morphisms
$g_s,h_s:Z_s\to X$ such that $g=g_s\circ p(Z)_s$, $h=h_s\circ p(Z)_s$
and $f\circ g_s=f\circ h_s$.
Since $f$ is a monomorphism of $D$, $g_s=h_s$ which implies
$g=h$.
\par C. The proof of 3) follows directly from 1) and 2).
\end{pf}

\begin{Prop}\label{XXX1.11} Suppose $f:X\to Y$ is a level morphism of pro-$C$
and $Z$ is an inverse system in $C$.
Let $\Sigma$ be the set of sequences $s$ in $I(X)$ such that
$(f_s)_\ast:Mor(Z,X_s)\to Mor(Z,Y_s)$ is a bijection.
If $\Sigma$
is cofinal in the set of all
sequences in $I(X)$, then $f_\ast:Mor(Z,X)\to Mor(Z,Y)$ is a bijection.
\end{Prop}
\begin{pf} Notice that it suffices to show that $f_\ast$
is a surjection.
Given $g:Z\to Y$ and $\alpha\in I(X)$
define $h_\alpha:Z\to X_\alpha$ as follows:

1. Find a sequence $s$ in $I(X)$ such that
$s(1)>\alpha$ and $(f_s)_\ast:Mor(Z,X_s)\to Mor(Z,Y_s)$ is a bijection.

2. Pick $k_s:Z\to X_s$ with $f_s\circ k_s=p(Y)_s\circ g$.

3. Define $h_\alpha$ as $p(X)^{s(1)}_\alpha\circ p(X_s)_1\circ k_s$.
\par Our first observation is that the above definition
does not depend on $s$. Indeed, if $t>s$,
then $f_s\circ (p(X)^t_s\circ k_t)=
p(Y)^t_s\circ f_t\circ k_t=
p(Y)^t_s\circ p(Y)_t\circ g= p(Y)_s\circ g=f_s\circ k_s$.
Therefore $p(X)^t_s\circ k_t=k_s$
and $p(X)^{s(1)}_\alpha\circ p(X_s)_1\circ k_s=
p(X)^{s(1)}_\alpha\circ p(X_s)_1\circ p(X)^t_s\circ k_t=
p(X)^{t(1)}_\alpha\circ p(X_t)_1\circ k_t$.
\par Using the first observation and given $\alpha<\beta$
one can find the same sequence $s$ to define both
$h_\alpha$ and $h_\beta$. Now,
$p(X)^\beta_\alpha\circ h_\beta=
p(X)^\beta_\alpha\circ p(X)^{s(1)}_\beta\circ p(X_s)_1\circ k_s=
p(X)^{s(1)}_\alpha\circ p(X_s)_1\circ k_s=h_\alpha$.
That means $\{h_\alpha\}_{\alpha\in I(X)}$
is a morphism $h$ from $Z$ to $X$.
\par Finally, for all $\alpha\in I(X)$,
$p(Y)_\alpha\circ (f\circ h)= f_\alpha\circ h_\alpha=
f_\alpha\circ p(X)^{s(1)}_\alpha\circ p(X_s)_1\circ k_s=
p(Y)^{s(1)}_\alpha\circ f_{s(1)}\circ p(X_s)_1\circ k_s=
p(Y)^{s(1)}_\alpha\circ p(Y_s)_1\circ (f_s\circ k_s)=
p(Y)^{s(1)}_\alpha\circ p(Y_s)_1\circ p(Y)_s\circ g=
p(Y)_\alpha\circ g$.
That means $f\circ h=g$.
\end{pf}

\medskip
\medskip

\section{Strong monomorphisms and strong epimorphisms}
\medskip
Unless stated otherwise, \(C\) is an arbitrary category in this section.
\par The following characterization of isomorphisms in pro-categories
is useful in introducing and understanding of the main concepts of this
section; strong monomorphisms and strong epimorphisms.

\begin{Prop}\label{XXX2.1}
Let \(f: X \rightarrow Y\) be a morphism in pro-\(C\).
\(f\) is an isomorphism if and only if for any \(P,Q \in Ob(C)\) and 
any commutative diagram

\begin{center}
\setlength{\unitlength}{1mm}
\begin{picture}(63,28)

\put(5.3,22){\(X\)}

\put(58.3,22){\(Y\)}

\put(5.3,0){\(P\)}

\put(58.3,0){\(Q\)}

\put(32,25){\(f\)}

\put(32,3){\(g\)}

\put(61,11){\(b\)}

\put(1,11){\(a\)}

\put(15,23){\vector(1,0){37}}

\put(15,1){\vector(1,0){37}}

\put(6.3,19){\vector(0,-1){13}}

\put(59.3,19){\vector(0,-1){13}}

\end{picture}
\end{center}

\noindent there exists  \(u:Y \rightarrow P\) such that
\( g \circ u= b\) and \( u \circ f=a\).
\end{Prop}
\begin{pf} If $f^{-1}$ exists, then clearly $u=a\circ f^{-1}$
satisfies the desired equalities.
\par Suppose morphism $u$ exists for any commutative diagram.
Without loss of generality, we may assume that 
\(f=\{f_{\alpha}:X_{\alpha} \rightarrow Y_{\alpha}\}_{\alpha \in 
I(X)}\) is a level morphism of pro-\(C\) from $X$
to $Y$.
Since
\begin{center}
\setlength{\unitlength}{1mm}
\begin{picture}(63,28)

\put(5.3,22){\(X\)}

\put(58.3,22){\(Y\)}

\put(5.3,0){\(X_\alpha\)}

\put(58.3,0){\(Y_\alpha\)}

\put(32,25){\(f\)}

\put(32,3){\(f_\alpha\)}

\put(61,11){\(p(Y)_\alpha\)}

\put(-6,11){\(p(X)_\alpha\)}

\put(15,23){\vector(1,0){37}}

\put(15,1){\vector(1,0){37}}

\put(6.3,19){\vector(0,-1){13}}

\put(59.3,19){\vector(0,-1){13}}

\end{picture}
\end{center}
\noindent
is commutative for any $\alpha\in A$,
there is $u_\alpha:Y\to X_\alpha$
such that $u_\alpha\circ f=p_\alpha$
and $f_\alpha\circ u_\alpha=p(Y)_\alpha$.
To prove that $f$ is an isomorphism
it suffices to show that $f^\ast:Mor(Y,P)\to Mor(X,P)$
is a bijection for each $P\in Ob(C)$
(see \ref{XXX1.5}).
Since any $g:X\to P$ factors
as $g=g'\circ p(X)_\alpha$ for some $\alpha\in I(X)$,
putting $h=g'\circ u_\alpha$ one gets
$h\circ f=g$, i.e. $f^\ast$ is a surjection.
If $g,h:Y\to P\in Ob(C)$ satisfy $f\circ g=f\circ h$,
then one can find representatives
$g',h':Y_\alpha\to P$ of $g$ and $h$
such that $g'\circ f_\alpha=h'\circ f_\alpha$.
Now $g=g'\circ q_\alpha=g'\circ f_\alpha\circ u_\alpha=
h'\circ f_\alpha\circ u_\alpha=h'\circ p(Y)_\alpha=h$,
i.e. $f^\ast$ is an injection.
\end{pf}

\begin{Def}\label{XXX2.2}
A morphism \(f: X \rightarrow Y\) in pro-\(C\) is called a {\bf 
strong monomorphism} ({\bf strong epimorphism}, respectively) if 
every commutative diagram
\begin{center}
\setlength{\unitlength}{1mm}
\begin{picture}(63,28)

\put(5.3,22){\(X\)}

\put(58.3,22){\(Y\)}

\put(5.3,0){\(P\)}

\put(58.3,0){\(Q\)}

\put(32,25){\(f\)}

\put(32,3){\(g\)}

\put(61,11){\(b\)}

\put(1,11){\(a\)}

\put(15,23){\vector(1,0){37}}

\put(15,1){\vector(1,0){37}}

\put(6.3,19){\vector(0,-1){13}}

\put(59.3,19){\vector(0,-1){13}}

\end{picture}
\end{center}
\noindent
  where \(P,Q \in Ob(C)\), admits a morphism \(u: Y\rightarrow P\) so that
\(u \circ f=a\)  (\(g \circ u = b\), respectively).
\end{Def}

\begin{Rem}\label{XXX2.2.5}
Notice that, provided $C$ has a terminal object $\ast$, the object
$Q$ in the above diagram is irrelevant (in the case of strong monomorphisms)
as it can always be
replaced by such $\ast$.
\end{Rem}

\begin{Rem}\label{XXX2.3}
If $X$ and $Y$ are objects of $C$, then $f:X\to Y$
is a strong monomorphism (strong epimorphism, respectively)
if and only if $f$ has a left inverse (a right inverse, respectively).
Simply put $g=f$, $a=id_X$, and $b=id_Y$ in the above diagram.
\end{Rem}

\begin{Rem}\label{XXX2.3.5}
Later on (see  \ref{XXX3.2.5} and  \ref{XXX3.12.5}) we will see examples of
strong monomorphisms (respectively, strong epimorphisms)
in the category pro-$Gr$ of pro-groups which do not have
a left (respectively, right) inverse.
\end{Rem}

\begin{Prop}\label{XXX2.4} Suppose that \(f=\{f_{\alpha}:X_{\alpha} 
\rightarrow Y_{\alpha}\}_{\alpha \in I(X)}\) is a level morphism of 
pro-\(C\) from $X$
to $Y$. For any commutative diagram
\begin{center}
\setlength{\unitlength}{1mm}
\begin{picture}(63,28)

\put(5.3,22){\(X\)}

\put(58.3,22){\(Y\)}

\put(5.3,0){\(P\)}

\put(58.3,0){\(Q\)}

\put(32,25){\(f\)}

\put(32,3){\(g\)}

\put(61,11){\(b\)}

\put(1,11){\(a\)}

\put(15,23){\vector(1,0){37}}

\put(15,1){\vector(1,0){37}}

\put(6.3,19){\vector(0,-1){13}}

\put(59.3,19){\vector(0,-1){13}}

\end{picture}
\end{center}
  \noindent we may find \(\alpha \in I(X)\) and representatives 
\(a_{\alpha}: X_{\alpha} \rightarrow P\) of $a$ and 
\(b_{\alpha}:Y_{\alpha} \rightarrow Q\)
of $b$ such that

\begin{center}
\setlength{\unitlength}{1mm}
\begin{picture}(63,28)

\put(5.3,22){\(X_{\alpha}\)}

\put(58.3,22){\(Y_{\alpha}\)}

\put(5.3,0){\(P\)}

\put(58.3,0){\(Q\)}

\put(32,25){\(f_{\alpha}\)}

\put(32,3){\(g\)}

\put(61,11){\(b_{\alpha}\)}

\put(1,11){\(a_{\alpha}\)}

\put(15,23){\vector(1,0){37}}

\put(15,1){\vector(1,0){37}}

\put(6.3,19){\vector(0,-1){13}}

\put(59.3,19){\vector(0,-1){13}}

\end{picture}
\end{center}
\noindent
is commutative.
\end{Prop}
\begin{pf} Choose representatives $u:X_\beta\to P$ of $a$
and $v:Y_\beta\to Q$ of $b$.
Since $g\circ u\circ p(X)_\beta=g\circ a=b\circ f=v\circ p(Y)_\beta\circ f=
v\circ f_\beta \circ p(X)_\beta$,
there is $\alpha>\beta$ such that
$g\circ u\circ p(X)_\beta^\alpha=
v\circ f_\beta \circ p(X)_\beta^\alpha$.
Put $a_\alpha=u\circ p(X)_\beta^\alpha$
and $b_\alpha=v\circ p(Y)_\beta^\alpha$.
\end{pf}

The following characterization of strong monomorphisms and
strong epimorphisms is especially useful. Its immediate consequence
is that both properties are preserved by functors
pro-$F$ if $F:C\to D$.
\begin{Prop}\label{XXX2.5}
Suppose that \(f=\{f_{\alpha}:X_{\alpha} \rightarrow 
Y_{\alpha}\}_{\alpha \in I(X)}\) is a level morphism of pro-\(C\). 
The following statements are equivalent:

a. \(f\) is a strong monomorphism (strong epimorphism, respectively).

b. For each
\(\alpha \in I(X)\) there is  a morphism
\(u_{\alpha}: Y \rightarrow X_{\alpha}\) such that
\(u_{\alpha} \circ f = p(X)_{\alpha}\)
(\( f_{\alpha} \circ u_{\alpha} = p(Y)_{\alpha}\), respectively).

c.  For each
\(\alpha \in I(X)\) there is \(\beta \in I(X)\), \(\beta > \alpha\) 
and a morphism
\(g_{\alpha,\beta}: Y_{\beta} \rightarrow X_{\alpha}\) such that
\(g_{\alpha, \beta} \circ f_{\beta} = p(X)_{\alpha}^{\beta}\)
(\( f_{\alpha} \circ g_{\alpha, \beta} = p(Y)_{\alpha}^{\beta}\), 
respectively).
\end{Prop}

\begin{pf} a)$\implies$b) follows from the definition
of strong monomorphisms (strong epimorphisms, respectively).
\par b)$\implies$c).  $u_{\alpha}$ has a representative
$v:Y_\gamma\to X_{\alpha}$ for some $\gamma>\alpha$.
Since \(u_{\alpha} \circ f = p(X)_{\alpha}\)
(\( f_{\alpha} \circ u_{\alpha} = p(Y)_{\alpha}\), respectively),
$v\circ f_\gamma$ and $p(X)_\alpha^\gamma$
($f_{\alpha} \circ v$ and $p(Y)_\alpha^\gamma$, respectively)
are representatives of the same morphism
from $X$ to $X_\alpha$ ($Y$ to $Y_\alpha$, respectively),
so there is $\beta>\gamma$ such that
$v\circ f_\gamma\circ p(X)^\beta_\gamma=p_\alpha^\gamma\circ p(X)^\beta_\gamma$
($f_{\alpha} \circ v\circ p(Y)^\beta_\gamma=p(Y)_\alpha^\gamma\circ 
p(Y)^\beta_\gamma$, respectively).
Put $g_{\alpha, \beta}=v\circ p(Y)^\beta_\gamma$.

\par c)$\implies$a).
Given a diagram

\begin{center}
\setlength{\unitlength}{1mm}
\begin{picture}(63,28)

\put(5.3,22){\(X\)}

\put(58.3,22){\(Y\)}

\put(5.3,0){\(P\)}

\put(58.3,0){\(Q\)}

\put(32,25){\(f\)}

\put(32,3){\(g\)}

\put(61,11){\(b\)}

\put(1,11){\(a\)}

\put(15,23){\vector(1,0){37}}

\put(15,1){\vector(1,0){37}}

\put(6.3,19){\vector(0,-1){13}}

\put(59.3,19){\vector(0,-1){13}}

\end{picture}
\end{center}

\noindent we may find \(\alpha \in I(X)\) and representatives 
\(a_{\alpha}: X_{\alpha} \rightarrow P\) and \(b_{\alpha}:Y_{\alpha} 
\rightarrow Q\) such that

\begin{center}
\setlength{\unitlength}{1mm}
\begin{picture}(63,28)

\put(5.3,22){\(X_{\alpha}\)}

\put(58.3,22){\(Y_{\alpha}\)}

\put(5.3,0){\(P\)}

\put(58.3,0){\(Q\)}

\put(32,25){\(f_{\alpha}\)}

\put(32,3){\(g\)}

\put(61,11){\(b_{\alpha}\)}

\put(1,11){\(a_{\alpha}\)}

\put(15,23){\vector(1,0){37}}

\put(15,1){\vector(1,0){37}}

\put(6.3,19){\vector(0,-1){13}}

\put(59.3,19){\vector(0,-1){13}}

\end{picture}
\end{center}
\noindent is commutative (see \ref{XXX2.4}). Let \(g_{\alpha,\beta}: 
Y_{\beta} \rightarrow X_{\alpha}\) be as in c)
and define \(u=a_{\alpha} \circ g_{\alpha,\beta} \circ p(Y)_{\beta}: 
Y \rightarrow P\).
Then, in the case of $f$ being a strong monomorphism,
\[ u \circ f=  a_{\alpha} \circ g_{\alpha,\beta} \circ p(Y)_{\beta} \circ f =
a_{\alpha} \circ g_{\alpha,\beta} \circ f_{\beta} \circ p(X)_{\beta}=\]

\[
=a_{\alpha} \circ p(X)_{\alpha}^{\beta} \circ p(X)_{\beta} = 
a_{\alpha} \circ p(X)_{\alpha}=a.
\]
Similarly, in the case of $f$ being a strong epimorphism,
\[ g \circ u= g \circ a_{\alpha} \circ g_{\alpha,\beta} \circ 
p(Y)_{\beta} =
b_{\alpha} \circ f_{\alpha} \circ g_{\alpha,\beta} \circ p(Y)_{\beta}=
\]

\[
=b_{\alpha} \circ p(Y)_{\alpha}^{\beta} \circ p(Y)_{\beta} = 
b_{\alpha} \circ p(Y)_{\alpha}=b.
\]
\end{pf}

\begin{Rem}\label{XXX2.6} In view of \ref{XXX2.5} one can relate 
strong epimorphisms and strong monomorphisms
to the following concepts previously discussed in literature:

1. {\bf Weak dominations} introduced by Dydak \cite{D1}
  (see also \cite{MS} p.186)
are precisely strong epimorphisms of pro-$H_0$.

2. Given a compact subset $X$ of the Hilbert cube $Q$ one considers
the system $N(X)$ of neighborhoods of $X$ in $Q$. It is an object of pro-$T$,
where $T$ is the category of topological spaces, and one has a natural
morphism $i:X\to N(X)$ of pro-$T$. Notice that $i$ is a strong monomorphism
of pro-$T$ if and only if $X$ is an ANR.

3. The above morphism $i:X\to N(X)$ can be interpreted
as a morphism of pro-$HT$, where $HT$ is the homotopy category
of topological spaces. Notice that $i$ is a strong epimorphism of
pro-$HT$ if and only if $X$ is {\bf internally movable} (see \cite{Bog}).
Indeed, $X$ is internally movable if for every neighborhood $U$ of $X$
in $Q$ there is a neighborhood $V$ of $X$ in $U$ and a map
$r:V\to X$ which is homotopic in $U$ to the inclusion $V\to U$.

4. {\bf Approximate ANR's in the sense of Clapp} \cite{Cla}
are introduced in a way related to strong monomorphisms.
Recall that $X\in AANR_C$ if for each $\epsilon>0$
there is a neighborhood $U$ of $X$ in $Q$ and a map
$r:U\to X$ such that $r|U$ is $\epsilon$-close to $id_X$.
Also, {\bf approximate ANR's in the sense of Noguchi} \cite{No}
  and AWNR's of \cite{Bog}, \cite{Ts} are defined in a way resembling
strong monomorphisms.

\end{Rem}

\begin{Prop}\label{XXX2.7}
Let \(f: X \rightarrow Y\) be a morphism in pro-\(C\). The following 
statements are equivalent:

1. \(f\) is a strong monomorphism.

2. $f^\ast:Mor(Y,P)\to Mor(X,P)$ is a surjection for each $P\in Ob(C)$.
\end{Prop}
\begin{pf} Without loss of generality we may
assume that \(f=\{f_{\alpha}:X_{\alpha} \rightarrow 
Y_{\alpha}\}_{\alpha \in I(X)}\) is a level morphism of pro-\(C\) 
from $X$
to $Y$.
\par 1)$\implies$2). Given $g:X\to P\in Ob(C)$ there is
a representative $v:X_\alpha\to P$ of $g$. By \ref{XXX2.5} there is
$u:Y\to X_\alpha$ such that $u\circ f=p(X)_\alpha$.
Put $h=v\circ u:Y\to P$. Now $h\circ f=v\circ u\circ f=v\circ p(X)_\alpha=g$.
\par 2)$\implies$1). Given $\alpha\in A$ there is $u:Y\to X_\alpha$
such that $u\circ f=p(X)_\alpha$. By \ref{XXX2.5}, $f$ is a strong 
monomorphism.
\end{pf}

\begin{Cor}\label{XXX2.8} If $f$ is a strong monomorphism (strong 
epimorphism, respectively)
of pro-$C$,
then $f$ is a monomorphism (epimorphism, respectively) of pro-$C$.
\end{Cor}
\begin{pf} Without loss of generality we may
assume that \(f=\{f_{\alpha}:X_{\alpha} \rightarrow 
Y_{\alpha}\}_{\alpha \in I(X)}\) is a level morphism of pro-\(C\) 
from $X$
to $Y$.
\par Suppose $f$ is a strong monomorphism
and $a,b:Z\to X$ are two morphisms of pro-$C$
such that $f\circ a=f\circ b$.
To show $a=b$ it suffices to prove $p(X)_\alpha\circ a=p(X)_\alpha\circ b$
for all $\alpha\in I(X)$.
Choose $u:Y\to X_\alpha$ such that $u\circ f=p(X)_\alpha$
(see \ref{XXX2.5}). Now $p(X)_\alpha\circ a=u\circ f\circ a=
u\circ f\circ b=p(X)_\alpha\circ b$.
\par Suppose $f$ is a strong epimorphism
and $a,b:Y\to Z$ are two morphisms of pro-$C$
such that $a\circ f=b\circ f$.
\par Special Case: $Z\in Ob(C)$. Choose
representatives $a',b':Y_\alpha\to Z$
of $a$ and $b$, respectively, such that
$a'\circ f_\alpha=b'\circ f_\alpha$.
By \ref{XXX2.5} there is $\beta>\alpha$ and $g:Y_\beta\to X_\alpha$
such that $f_\alpha\circ g=p(Y)^\beta_\alpha$.
Therefore $a'\circ p(Y)^\beta_\alpha= a'\circ f_\alpha\circ g=
b'\circ f_\alpha\circ g=b'\circ p(Y)^\beta_\alpha$
which proves $a=b$.
\par General Case.
To show $a=b$ we need $p(Z)_i\circ a=p(Z)_i\circ b$ for
all $i\in I(Z)$ which follows from Special Case.
\end{pf}

\begin{Rem}\label{XXX2.8.5} Notice that there are
monomorphisms (respectively, epimorphisms) of pro-$Gr$ 
that are not strong monomorphisms (respectively,
strong epimorphisms). Easy examples are the inclusion
$Z\to Q$ from integers to rational numbers and the projection
$Z\to Z/2$ from integers to integers modulo $2$.
\end{Rem}

\begin{Cor}\label{XXX2.9} If $g\circ f$ is a strong monomorphism 
(strong epimorphism, respectively),
then $f$ is a strong monomorphism ($g$ is a strong epimorphism, respectively).
\end{Cor}
\begin{pf} Assume $f:X\to Y$ and $g:Y\to Z$.
\par
Suppose $g\circ f$ is a strong monomorphism
and $a:X\to P\in Ob(C)$. By \ref{XXX2.7} there is $b:Z\to P$
such that $a=b\circ (g\circ f)$.
Now, $c=b\circ g$ satisfies $a=c\circ f$
and \ref{XXX2.7} says that $f$ is a strong monomorphism.
\par Suppose $g\circ f$ is a strong epimorphism
and
\begin{center}
\setlength{\unitlength}{1mm}
\begin{picture}(63,28)

\put(5.3,22){\(Y\)}

\put(58.3,22){\(Z\)}

\put(5.3,0){\(P\)}

\put(58.3,0){\(Q\)}

\put(32,25){\(g\)}

\put(32,3){\(h\)}

\put(61,11){\(b\)}

\put(1,11){\(a\)}

\put(15,23){\vector(1,0){37}}

\put(15,1){\vector(1,0){37}}

\put(6.3,19){\vector(0,-1){13}}

\put(59.3,19){\vector(0,-1){13}}

\end{picture}
\end{center}

\noindent
is a commutative diagram in pro-$C$ with $P,Q\in Ob(C)$.
Since
\begin{center}
\setlength{\unitlength}{1mm}
\begin{picture}(63,28)

\put(5.3,22){\(X\)}

\put(58.3,22){\(Z\)}

\put(5.3,0){\(P\)}

\put(58.3,0){\(Q\)}

\put(28,25){\(g\circ f\)}

\put(32,3){\(h\)}

\put(61,11){\(b\)}

\put(-3,11){\(a\circ f\)}

\put(15,23){\vector(1,0){37}}

\put(15,1){\vector(1,0){37}}

\put(6.3,19){\vector(0,-1){13}}

\put(59.3,19){\vector(0,-1){13}}

\end{picture}
\end{center}
\noindent
is commutative, there is $u:Z\to P$
with $h\circ u=b$
which proves that $g$ is a strong epimorphism.
\end{pf}

The following is the main property of strong monomorphisms.

\begin{Thm}\label{XXX2.8a}
Suppose
\begin{center}
\setlength{\unitlength}{1mm}
\begin{picture}(63,28)

\put(5.3,22){\(Z\)}

\put(58.3,22){\(T\)}

\put(5.3,0){\(X\)}

\put(58.3,0){\(Y\)}

\put(32,25){\(f'\)}

\put(32,3){\(f\)}

\put(61,11){\(g\)}

\put(1,11){\(g'\)}

\put(15,23){\vector(1,0){37}}

\put(15,1){\vector(1,0){37}}

\put(6.3,19){\vector(0,-1){13}}

\put(59.3,19){\vector(0,-1){13}}

\end{picture}
\end{center}
\noindent
is a commutative diagram in pro-$C$.
If $f'$ is an epimorphism and $f$ is a strong monomorphism,
then there is a unique filler $u:T\to X$, i.e.
a morphism $u$ such that $g'=u\circ f'$ and $f\circ u=g$.
\end{Thm}
\begin{pf} Since $f'$ is an epimorphism, it suffices
to prove existence of $u$.
Also, it suffices to prove that $g'=u\circ f'$.
Indeed, $g'=u\circ f'$ implies $(f\circ u) \circ f'=f\circ g'=g\circ f'$,
so $f\circ u=g$ as $f'$ is an epimorphism.
\par Given $\alpha\in I(X)$ there is $r(\alpha):Y\to X_\alpha$
such that $r(\alpha)\circ f=p(X)_\alpha$.
Put $u_\alpha=r(\alpha)\circ g$.
If $\beta>\alpha$, then $p(X)^\beta_\alpha\circ u_\beta\circ f'=
  p(X)^\beta_\alpha\circ r(\beta)\circ g\circ f'=
p(X)^\beta_\alpha\circ r(\beta)\circ f\circ g'=
p(X)^\beta_\alpha\circ p(X)_\beta\circ  g'=
p(X)_\alpha\circ  g'$.
Similarly,
$u_\alpha\circ f'=
p(X)_\alpha\circ  g'$. Since $f'$ is an epimorphism,
$u_\alpha=p(X)^\beta_\alpha\circ u_\beta$ which means that 
$\{u_\alpha\}_{\alpha\in I(X)}$
is a morphism from $T$ to $X$.
Also, $u_\alpha\circ f'=
p(X)_\alpha\circ  g'$ for all $\alpha\in I(X)$ means $u\circ f'=g'$.
\end{pf}

\begin{Cor} \label{XXX2.10}
Let \(f: X \rightarrow Y\) be a morphism in pro-\(C\). The following 
statements are equivalent:

1. \(f\) is an isomorphism.

2. \(f\) is a strong monomorphism and an epimorphism.
\end{Cor}
\begin{pf}
Obviously, only $2)\implies 1)$ is of interest.
Apply \ref{XXX2.8a}
to the diagram
\begin{center}
\setlength{\unitlength}{1mm}
\begin{picture}(63,28)

\put(5.3,22){\(X\)}

\put(58.3,22){\(Y\)}

\put(5.3,0){\(X\)}

\put(58.3,0){\(Y\)}

\put(32,25){\(f\)}

\put(32,3){\(f\)}

\put(61,11){\(id_Y\)}

\put(-1,11){\(id_X\)}

\put(15,23){\vector(1,0){37}}

\put(15,1){\vector(1,0){37}}

\put(6.3,19){\vector(0,-1){13}}

\put(59.3,19){\vector(0,-1){13}}

\end{picture}
\end{center}
\end{pf}

The following is the main property of strong epimorphisms.

\begin{Thm}\label{XXX2.12a} Let $C$ be a category
with direct sums.
Suppose
\begin{center}
\setlength{\unitlength}{1mm}
\begin{picture}(63,28)

\put(5.3,22){\(Z\)}

\put(58.3,22){\(T\)}

\put(5.3,0){\(X\)}

\put(58.3,0){\(Y\)}

\put(32,25){\(f'\)}

\put(32,3){\(f\)}

\put(61,11){\(g\)}

\put(1,11){\(g'\)}

\put(15,23){\vector(1,0){37}}

\put(15,1){\vector(1,0){37}}

\put(6.3,19){\vector(0,-1){13}}

\put(59.3,19){\vector(0,-1){13}}

\end{picture}
\end{center}
\noindent
is a commutative diagram in pro-$C$.
If $f'$ is a strong epimorphism and $f$ is a monomorphism,
then there is a unique filler $u:T\to X$, i.e.
a morphism $u$ such that $g'=u\circ f'$ and $f\circ u=g$.
\end{Thm}
\begin{pf} Since $f$ is a monomorphism, it suffices
to prove existence of $u$.
Also, it suffices to prove that $g=f\circ u$.
Indeed, $g=f\circ u$ implies $f\circ (u \circ f')=g\circ f'=f\circ g'$,
so $u\circ f'=g'$ as $f$ is a monomorphism.
\par Assume $f:X\to Y$ is a level morphism
such that $I(X)$ is cofinite.
Let $n(\alpha)$ be the number of predecessors
of $\alpha\in I(X)$. By induction on $n(\alpha)$
one can find an increasing function $e:I(X)\to I(X)$
such that for any two morphisms
$a,b:Z\to X_{e(\alpha)}$ the equality
$f_{e(\alpha)}\circ a=f_{e(\alpha)}\circ b$
implies $p(X)^{e(\alpha)}_{\alpha}\circ a=p(X)^{e(\alpha)}_{\alpha}\circ b$
(see \ref{XXX1.mono}).
\par Since the diagram
\begin{center}
\setlength{\unitlength}{1mm}
\begin{picture}(63,50)

\put(5.3,44){\(Z\)}

\put(58.3,44){\(T\)}

\put(5.3,22){\(X\)}

\put(58.3,22){\(Y\)}

\put(5.3,0){\(X_\alpha\)}

\put(58.3,0){\(Y_\alpha\)}

\put(32,47){\(f'\)}

\put(32,25){\(f\)}

\put(32,3){\(f_\alpha\)}

\put(61,33){\(g\)}

\put(-1,33){\(g'\)}

\put(61,11){\(p(Y)_\alpha\)}

\put(-5,11){\(p(X)_\alpha\)}

\put(15,45){\vector(1,0){37}}

\put(15,23){\vector(1,0){37}}

\put(15,1){\vector(1,0){37}}

\put(6.3,41){\vector(0,-1){13}}

\put(59.3,41){\vector(0,-1){13}}

\put(6.3,19){\vector(0,-1){13}}

\put(59.3,19){\vector(0,-1){13}}

\end{picture}
\end{center}
\noindent
is commutative, one has $h_\alpha:T\to X_\alpha$
so that $f_\alpha\circ h_\alpha= p(Y)_\alpha\circ g$.
Define $u_\alpha=p(X)^{e(\alpha)}_\alpha\circ h_{e(\alpha)}$.
Suppose $\beta>\alpha$. Notice that
  $f_{e(\alpha)}\circ p(X)^{e(\beta)}_{e(\alpha)}\circ h_{e(\beta)}=
p(Y)^{e(\beta)}_{e(\alpha)}\circ f_{e(\beta)}\circ h_{e(\beta)}=
p(Y)^{e(\beta)}_{e(\alpha)}\circ p(Y)_{e(\beta)}\circ g=
p(Y)_{e(\alpha)}\circ g$.
Also, $f_{e(\alpha)}\circ h_{e(\alpha)}=
p(Y)_{e(\alpha)}\circ g$,
so $p(X)^{e(\alpha)}_{\alpha}\circ h_{e(\alpha)}
=p(X)^{e(\alpha)}_{\alpha}\circ p(X)^{e(\beta)}_{e(\alpha)}\circ h_{e(\beta)}$,
i.e. $u_\alpha=p(X)^\beta_\alpha\circ u_\beta$.
That means $u=\{u_\alpha\}_{\alpha\in I(X)}$
is a morphism from $T$ to $X$.
Since $f_\alpha\circ u_\alpha=
f_\alpha\circ p(X)^{e(\alpha)}_{\alpha}\circ h_{e(\alpha)}=
p(Y)^{e(\alpha)}_{\alpha}\circ f_{e(\alpha)}\circ h_{e(\alpha)}
=p(Y)^{e(\alpha)}_{\alpha}\circ p(Y)_{e(\alpha)}\circ g=
p(Y)_{\alpha}\circ g$ for each $\alpha\in I(X)$,
we have $f\circ u=g$.
\end{pf}


\begin{Cor} \label{XXX2.12}
Let \(f: X \rightarrow Y\) be a morphism in pro-\(C\). If \(C\) has 
direct sums, then the following statements are equivalent:

1. \(f\) is an isomorphism.

2. \(f\) is a strong epimorphism and a monomorphism.
\end{Cor}
\begin{pf} Obviously, only $2)\implies 1)$ is of interest.
Apply \ref{XXX2.12a}
to the diagram
\begin{center}
\setlength{\unitlength}{1mm}
\begin{picture}(63,28)

\put(5.3,22){\(X\)}

\put(58.3,22){\(Y\)}

\put(5.3,0){\(X\)}

\put(58.3,0){\(Y\)}

\put(32,25){\(f\)}

\put(32,3){\(f\)}

\put(61,11){\(id_Y\)}

\put(-1,11){\(id_X\)}

\put(15,23){\vector(1,0){37}}

\put(15,1){\vector(1,0){37}}

\put(6.3,19){\vector(0,-1){13}}

\put(59.3,19){\vector(0,-1){13}}

\end{picture}
\end{center}
\end{pf}

\medskip
\medskip

\section{Movability}

\medskip

In this section we introduce a new variant of movability
and we discuss the connection of movability to various
classes of morphisms of pro-categories.

\begin{Def}\label{XXX3.1}
  \(X \in Ob(pro-C)\) is {\bf uniformly movable}
if the morphism $\{p(X)_\alpha:X\to X_\alpha\}_{\alpha\in I(X)}$
from $X$ to $Sub_2(X)$
is a strong epimorphism of pro-(pro-$C$).
In view of \ref{XXX2.5} it means that for each $\alpha\in I(X)$
there is $\beta>\alpha$ and $r:X_\beta\to X$
such that $r_\alpha=p(X)^\beta_\alpha$,
which is the classical definition of uniform movability
(see \cite{MS}, p.160).
\end{Def}

Here is the connection between uniform movability and strong
epimorphisms.
\begin{Prop}\label{XXX3.2}

a. If there is \(P \in Ob(C)\) and a strong epimorphism \(f:P \rightarrow X\),
then  \(X\) is uniformly movable.

b. If $C$ is a category with inverse limits
and $X$ is uniformly movable, then the projection
$p:\underset{\leftarrow}{\text{lim}}(X)\to X$ is a strong epimorphism.

c. If \(C\) is a category with direct sums and \(X \in Ob(pro-C)\) is 
uniformly movable,
then there is \(P \in Ob(C)\) and a strong epimorphism \(P \rightarrow X\).
\end{Prop}
\begin{pf}
a. Suppose that \(f: P \rightarrow X\) is a strong epimorphism and 
\(\alpha \in I(X)\).
There is \(\beta> \alpha\) and \(g: X_{\beta} \rightarrow P\) so that 
\(f_{\alpha} \circ g =
p(X)_{\alpha}^{\beta}\). Now, \( p(X)_{\alpha} \circ f \circ g= 
f_{\alpha} \circ g=p(X)_{\alpha}^{\beta}\) which means that \(X\) is 
uniformly movable.
\par  In b) and c) assume that for each \(\alpha \in I(X)\) there is 
\(\beta(\alpha) > \alpha\) and
\(g_{\alpha}: X_{\beta(\alpha)} \rightarrow X\) so that
\(p(X)_{\alpha} \circ g_{\alpha} =p(X)_{\alpha}^{\beta(\alpha)}\).
\par b. $g_{\alpha}$ factors through $\underset{\leftarrow}{\text{lim}}(X)$
so that there is $h_{\alpha}:X_\beta\to \underset{\leftarrow}{\text{lim}}(X)$
with $g_\alpha=p\circ h_\alpha$. Now
$(p(X)_\alpha\circ p)\circ h_\alpha=p(X)^\beta_\alpha$
which proves that the level morphism
$\{p(X)_\alpha\circ p\}_{\alpha\in I(X)}$
is a strong epimorphism (see \ref{XXX2.5}).
That morphism is exactly $p:\underset{\leftarrow}{\text{lim}}(X)\to X$.
\par c. Let \(P=\bigoplus\limits_{\alpha \in I(X)} 
X_{\beta(\alpha)}\) and let \(f:P \rightarrow X\) be induced by 
\(g_{\alpha}\), \( \alpha \in I(X)\).
Given \(\alpha \in I(X)\) we have \(i_{\alpha} : X_{\beta(\alpha)} 
\rightarrow P\) so that \(f_{\alpha} \circ i_{\alpha}= p(X)_{\alpha} 
\circ f \circ i_{\alpha} = p(X)_{\alpha} \circ g_{\alpha}= 
p(X)_{\alpha}^{\beta(\alpha)}\). That means \(f\) is a strong 
epimorphism.
\end{pf}

\begin{Rem}\label{XXX3.2.5} Consider a uniformly movable
pro-group $G$ which is not stable. By \ref{XXX3.2} there
is a strong epimorphism $f:P\to G$ from a group $P$. That epimorphism
cannot have a right inverse as $G$ is not stable.
\end{Rem}

In \cite{MS} (Theorem 3 on p.162)
it is shown that if $Y$ is uniformly movable and
$\underset{\leftarrow}{\text{lim}}(f):\underset{\leftarrow}{\text{lim}}(X)\to 
\underset{\leftarrow}{\text{lim}}(Y)$
is an epimorphism of $C$, then $f$ is an epimorphism of
pro-$C$. We derive that result in part a) below.
\begin{Cor}\label{XXX3.3}
Suppose \(C\) is a category with inverse limits and \(f: X \rightarrow Y\) is
morphism of pro-\(C\) such that $Y$ is uniformly movable.

a. If $\underset{\leftarrow}{\text{lim}}(f)$ is
an epimorphism of $C$, then $f$ is an epimorphism of pro-$C$.

b. If $\underset{\leftarrow}{\text{lim}}(f)$ 
has a right inverse in $C$, then $f$ is a strong epimorphism of pro-$C$.

c. If $C$ is a category with direct sums, $f$ is a monomorphism
of pro-$C$, and
$\underset{\leftarrow}{\text{lim}}(f)$ is
an isomorphism of $C$, then $f$ is an isomorphism of pro-$C$.
\end{Cor}

\begin{pf} a. (respectively, b.)
By \ref{XXX3.2} the projection morphism 
$\underset{\leftarrow}{\text{lim}}(Y)\to Y$
is a strong epimorphism. Therefore the composition
$\underset{\leftarrow}{\text{lim}}(X) \rightarrow 
\underset{\leftarrow}{\text{lim}}(Y)\to Y$
is an epimorphism (respectively, a strong epimorphism by 
\ref{XXX2.9}) of pro-$C$.
That composition equals
$\underset{\leftarrow}{\text{lim}}(X) \rightarrow X\overset{f}\rightarrow Y$,
so $f$ is an epimorphism (respectively, a strong epimorphism).
\par c. It follows from b) and \ref{XXX2.12}.
\end{pf}

Notice that, for a category $C$ which does not have
inverse limits, the analog of $\underset{\leftarrow}{\text{lim}}(f)$
being an isomorphism (see \ref{XXX3.3}) is that
$f_\ast:Mor(P,X)\to Mor(P,Y)$ is a bijection for all $P\in Ob(C)$.
Our next results should be viewed in that context.

\begin{Cor}\label{XXX3.4}
Suppose \(C\) is a category with direct sums or inverse limits
and \(Y \in Ob(pro-C)\) is uniformly movable.
A monomorphism $f:X\to Y$ of $pro-C$ is an isomorphism if and only if
$f_\ast:Mor(P,X)\to Mor(P,Y)$ is a surjection for each $P\in Ob(C)$.
\end{Cor}
\begin{pf} Find a strong epimorphism $g:P\to Y$ such that $P\in Ob(C)$.
Factor $g$ as $f\circ h$ for some $h:P\to X$.
By \ref{XXX2.9}, $f$ is a strong epimorphism, so \ref{XXX2.12} implies that
$f$ is an isomorphism.
\end{pf}

\begin{Thm}\label{XXX3.15a} Suppose $C$ is a category with direct sums
and $f:X\to Y$ is a monomorphism of pro-$C$.
If $f_\ast:Mor(P,X)\to Mor(P,Y)$ is a surjection for all $P\in Ob(C)$,
then $f_\ast:Mor(T,X)\to Mor(T,Y)$ is a bijection for all
uniformly movable objects $T$ of pro-$C$.
\end{Thm}
\begin{pf} Pick a strong epimorphism $f':P\to T$
such that $P\in Ob(C)$ (see \ref{XXX3.2}).
Suppose $g:T\to Y$ is a morphism and find
$g':P\to X$ such that $f\circ g'=g\circ f'$.
That means the diagram
\begin{center}
\setlength{\unitlength}{1mm}
\begin{picture}(63,28)

\put(5.3,22){\(P\)}

\put(58.3,22){\(T\)}

\put(5.3,0){\(X\)}

\put(58.3,0){\(Y\)}

\put(32,25){\(f'\)}

\put(32,3){\(f\)}

\put(61,11){\(g\)}

\put(1,11){\(g'\)}

\put(15,23){\vector(1,0){37}}

\put(15,1){\vector(1,0){37}}

\put(6.3,19){\vector(0,-1){13}}

\put(59.3,19){\vector(0,-1){13}}

\end{picture}
\end{center}
\noindent
is commutative, so there is $u:T\to X$
such that $f\circ u=g$ by \ref{XXX2.12a}.
\end{pf}

\begin{Prop}\label{XXX3.5}
Suppose \(C\) is a category and \(X \in Ob(pro-C)\) is uniformly movable.
$X$ is dominated by an object of $C$
if there is a monomorphism $f:X\to P$, where $P\in Ob(C)$.
\end{Prop}
\begin{pf}
Pick $g:X_\alpha\to P$ for some $\alpha\in I(X)$, so that
$f=g\circ p(X)_\alpha$. Choose $h:X_\beta\to X$
with $p(X)_\alpha\circ h=p(X)^\beta_\alpha$.
Let $u=h\circ p(X)_\beta:X\to X$. We plan to show that $u=id_X$
which means that $X$ is dominated by $X_\beta$.
It suffices to show that $f\circ u=f\circ id_X$
as $f$ is a monomorphism. Indeed,
$f\circ u= f\circ h\circ p(X)_\beta=g\circ p(X)_\alpha \circ h\circ p(X)_\beta=
g\circ p(X)^\beta_\alpha \circ p(X)_\beta=g\circ p(X)_\alpha=f$.
\end{pf}

\begin{Rem}\label{XXX3.5.5} In \ref{XXX3.11} we will see
a pro-group $G$ admitting a monomorphism to a group
such that $G$ is not stable.
Thus, the assumption of $G$ being uniformly movable is essential.
\end{Rem}


\begin{Prop}\label{XXX3.14a} Suppose $f:X\to Y$ is a morphism of pro-$C$
  such that
$f_\ast:Mor(P,X)\to Mor(P,Y)$ is a surjection
for all $P\in Ob(C)$. If $Y$ is uniformly movable, then $f$ is a 
strong epimorphism.
\end{Prop}
\begin{pf} Suppose $f$ is a level morphism.
Given $\alpha\in I(X)$ find $\beta>\alpha$
and $r:Y_\beta\to Y$
satisfying $r_\alpha=p(Y)^\beta_\alpha$.
Lift $r$ to $X$, i.e. find $h:Y_\beta\to X$
with $f\circ h=r$. Notice that
$f_\alpha\circ h_\alpha=r_\alpha =p(Y)^\beta_\alpha$
which means $f$ is a strong epimorphism by \ref{XXX2.5}.
\end{pf}

\begin{Rem}\label{XXX3.14.5} The pro-group $G$ from \ref{XXX3.11}
has the property that the inclusion $0\to G$ from the trivial group
induces epimorphisms $f_\ast:Mor(P,0)\to Mor(P,G)$ for all groups
$P$ but $0\to G$ is not a strong epimorphism.
Thus, the assumption of $G$ being uniformly movable is essential.
\end{Rem}

\begin{Def}\label{XXX3.6} \(X \in Ob(pro-C)\) is {\bf sequentially movable}
if the morphism $Sub_{\omega_0}(X)\to Sub_2(X)$ 
is a strong epimorphism of pro-(pro-$C$).
Alternatively, for any increasing sequence $s$ in $I(X)$
there is $\beta>s(1)$ and a morphism $r:X_\beta\to X_s$
such that $r_1=p(X)^\beta_{s(1)}$.
  \end{Def}

\begin{Rem}\label{XXX3.7} Notice that if $X$ is sequentially movable, then
  for any increasing sequence $s$ in $I(X)$ and any $k\ge 1$
there is $\beta>s(k)$ and a morphism $r:X_\beta\to X_s$
such that $r_k=p(X)^\beta_{s(k)}$.
\end{Rem}

\begin{Prop}\label{XXX3.8} If $X$ is a movable
object of pro-$C$, then it is sequentially movable.
\end{Prop}
\begin{pf} Clearly, if $X$ is uniformly movable,
then it is sequentially movable.
Notice that for every sequence $s$ in $I(X)$
there is a sequence $t>s$ such that $X_t$ is movable,
hence uniformly movable (see \cite{Sp} or \cite{MS}, Theorem 4
on p.163). Thus, there is $\beta=t(k)$ for some
$k>1$ and $u:X_\beta\to X_t$ with
$u_1=p(X)^\beta_{t(1)}$. Put $r=p(X)^t_s\circ u:X_\beta\to X_s$
to get $r_1=p(X)^\beta_{s(1)}$.
\end{pf}

\begin{Prop}\label{XXX3.9} If $X$
has the property that each morphism $p(X)_s$
from $X$ to its subtower $X_s$ factors through a sequentially movable
object of pro-$C$, then $X$
is sequentially movable.
\end{Prop}
\begin{pf} We need to factor through sequentially movable
objects whose index set is cofinite.
\par Claim. Every sequentially movable object $Z$
is isomorphic to a sequentially movable object $Z'$
such that $I(Z')$ is cofinite.
\par Proof of Claim. We will employ the standard
reindexing trick as in the proof of Theorem 3 on p.12 in \cite{MS}.
Define $I(Z')$ to be the set of
all finite subsets $\sigma$ of $I(Z)$ which have a maximum $\max(\sigma)$
and declare
$\sigma\leq \tau$ if $\sigma\subset \tau$.
$Z'_\sigma$ is defined as $Z_{\max(\sigma)}$
and $p(Z')^\tau_\sigma:=p(Z)_{\max(\sigma)}
^{\max(\tau)}$.
Given an increasing sequence $s$ in $I(Z')$
we define $t(n)=\max(s(n))$.
There is $\gamma>t(1)$ and $u:Z_\gamma\to Z_t$
such that $u_1=p(Z)^\gamma_{t(1)}$.
Setting $\beta=s(1)\cup\{\gamma\}$
and interpreting $u$ as $r:Z'_\beta\to Z'_s$
one gets $r_1=p(Z')^\beta_{s(1)}$.
Thus $Z'$ is sequentially movable.
Notice that projections from $Z'$ to $Z_\alpha$,
where $\alpha$ is interpreted as a one-point set,
form a morphism from $Z'$ to $Z$ which is an isomorphism.
\par Suppose $s$ is an increasing sequence in $I(X)$.
Factor $p(X)_s:X\to X_s$ as $h\circ g$,
where $g:X\to Z$, $h:Z\to X_s$, $Z$ is sequentially
movable, and $I(Z)$ is cofinite.
Find a sequence $t$ in $I(Z)$ and a level morphism
$f:Z_t\to X_s$ such that $h=f\circ p(Z)_t$.
For each $\alpha\in I(Z)$ let $n(\alpha)$ be the number
of predecessors of $\alpha$. By induction on $n(\alpha)$
we can construct an increasing function $i:I(Z)\to I(X)$
and representatives $g_\alpha:X_{i(\alpha)}\to Z_\alpha$
of $p(Z)_\alpha\circ g$ such that
$\{g_\alpha\}_{\alpha\in I(Z)}$ induces a level morphism,
i.e. $p(Z)^\beta_\alpha\circ g_\beta=g_\alpha\circ p(X)^{i(\beta)}_{i(\alpha)}$
for all $\alpha<\beta$.
Pick $\gamma\in I(Z)$, $\gamma>t(1)$, and $u:Z_\gamma\to Z_t$
such that $u_1=p(Z)^\gamma_{t(1)}$. Set $\beta=i(\gamma)$
and $r= f\circ u\circ g_\gamma$.
Now $p(X_s)_1\circ r=p(X_s)_1\circ f\circ u\circ g_\gamma=
f_1\circ p(Z_t)_1\circ u\circ g_\gamma=
f_1\circ p(Z)^\gamma_{t(1)}\circ g_\gamma=
f_1\circ g_{t(1)}\circ p(X)^\beta_{t(1)}=
p(X)^{i(t(1)}_{s(1)}\circ p(X)^\beta_{t(1)}=
p(X)^\beta_{s(1)}$ and $X$ is sequentially movable.
\end{pf}

\begin{Cor}\label{XXX3.10} If $X$ is dominated by a sequentially movable
object of pro-$C$, then it is sequentially movable.
\end{Cor}

Let us show that sequential movability is a more general concept
than movability.
\begin{Prop}\label{XXX3.11} There is a sequentially movable
pro-group which is not movable.
\end{Prop}
\begin{pf} Let $I(X)$ be the set of all
ordinals smaller than the first uncountable ordinal.
For each $\alpha\in I(X)$ let
$X_\alpha=\bigoplus\limits_{\gamma\in I(X)} G^\alpha_\gamma$,
where $G^\alpha_\gamma=0$ if $\gamma<\alpha$
and $G^\alpha_\gamma$ is the group of natural numbers $Z$ otherwise.
$p(X)^\beta_\alpha:X_\beta\to X_\alpha$ is the natural inclusion.
Any increasing sequence $s$ in $I(X)$ has an upper bound
$\beta$. Therefore any morphism from $X$
to a tower factors through a group. By \ref{XXX3.9} that proves the 
sequential movability
of $X$.
\par If $X$ were movable, then for each $\alpha$
there would be $\beta>\alpha$ with $im(p(X)^\beta_\alpha)$
contained in each $im(p(X)^\gamma_\alpha)$,
$\gamma>\beta$. However, that implies
$im(p(X)^\beta_\alpha)=0$, a contradiction.
Notice that \ref{XXX3.16} generalizes the above argument.
\end{pf}

\begin{Rem}\label{XXX3.12} The same argument as in
\ref{XXX3.11} shows that the system
of subtowers $Sub_{\omega_0}(X)$ of $X$ is always sequentially
movable. That is because every sequence in $Sub_{\omega_0}(X)$
has an upper bound.
\end{Rem}

\begin{Rem}\label{XXX3.12.5} Consider the inclusion $X\to X_0$ as in
\ref{XXX3.11}. It is a strong monomorphism of pro-$Gr$
which does not have a left inverse as $X$ is not stable.
\end{Rem}

\begin{Prop}\label{XXX3.13} If $f:X\to Y$ is a strong epimorphism
and $X$ is sequentially movable, then $Y$
is sequentially movable.
\end{Prop}
\begin{pf}
Without loss of generality assume that
$f$ is a level morphism. Given a sequence $s$ in $I(Y)$
find $\gamma\in I(X)$, $\gamma>s(1)$,
and $u:X_\gamma\to X_s$ such that $u_1=p(X)^\gamma_{s(1)}$.
Find $\beta>\gamma$ and $g:Y_\beta\to X_\gamma$
such that $f_\gamma\circ g=p(Y)^\beta_\gamma$.
Set $r=f_s\circ u\circ g:Y_\beta\to Y_s$.
Notice that $r_1=f_1\circ u_1\circ g=
f_1\circ p(X)^\gamma_{s(1)}\circ g=
p(Y)^\gamma_{s(1)}\circ f_\gamma\circ g=
p(Y)^\gamma_{s(1)}\circ p(Y)^\beta_\gamma=
p(Y)^\beta_{s(1)}$.
\end{pf}

\begin{Prop}\label{XXX3.14} Suppose $f:X\to Y$ is a level morphism of pro-$C$
such that $Y$ is sequentially movable.
Let $\Sigma$ be the set of sequences $s$ in $I(X)$ such that
$(f_s)_\ast:Mor(P,X_s)\to Mor(P,Y_s)$ is a surjection
for all $P\in Ob(C)$. If $\Sigma$ is cofinal in the set of all
sequences in $I(X)$, then $f$ is a strong epimorphism.
\end{Prop}
\begin{pf}
Given $\alpha\in I(X)$ find
a sequence $s$ in $I(X)$ such that $s(1)>\alpha$
and $(f_s)_\ast:Mor(P,X_s)\to Mor(P,Y_s)$ is a surjection
for all $P\in Ob(C)$. Pick $\beta>\alpha$ and $r:Y_\beta\to Y_s$
satisfying $r_1=p(Y)^\beta_{s(1)}$.
Lift $r$ to $X_s$, i.e. find $r':Y_\beta\to X_s$
with $f_s\circ r'=r$. Set $g:Y_\beta\to X_\alpha$ to be
$p(X)^{s(1)}_\alpha\circ p(X_s)_1\circ r'$
and notice that
$f_\alpha\circ g=f_\alpha \circ p(X)^{s(1)}_\alpha\circ p(X_s)_1\circ r'=
p(Y)^{s(1)}_\alpha\circ p(Y_s)_1\circ f_s\circ r'=
p(Y)^{s(1)}_\alpha\circ p(Y_s)_1\circ r=p(Y)^\beta_\alpha$.
\end{pf}

\begin{Thm}\label{XXX3.15} Suppose $C$ is a category with direct sums
and $f:X\to Y$ is a monomorphism of $tow(C)$.
If $f_\ast:Mor(P,X)\to Mor(P,Y)$ is a surjection for all $P\in Ob(C)$,
then $f_\ast:Mor(Z,X)\to Mor(Z,Y)$ is a bijection for all
sequentially movable objects $Z$ of pro-$C$.
\end{Thm}
\begin{pf} Special Case: $f$ is a level morphism induced by
$\{f_n:X_n\to Y_n\}_{n\in N}$ such that for any two morphisms
$a,b:P\to X_{n+1}$ of $C$ the equality $f_{n+1}\circ a=f_{n+1}\circ b$
implies $p(X)^{n+1}_n\circ a=p(X)^{n+1}_n\circ b$.
\par Suppose $g:Z\to Y$ is a morphism of pro-$C$
and $Z$ is sequentially movable.
By \ref{XXX1.8} there is an increasing sequence $s$ in $I(Z)$ and a level
morphism $h:Z_s\to Y$ satisfying $g=h\circ p(Z)_s$.
By induction on $n$ we can construct an increasing sequence
$t$ in $I(Z)$, $t>s$,
and morphisms $r(n):Z_{t(n)}\to Z_s$, $n\ge 1$,
so that $r(n)_n=p(Z)^{t(n)}_{s(n)}$.
Lift each $h\circ r(n)$ to $X$, i.e. find
$h(n):Z_{t(n)}\to X$ so that $f\circ h(n)=h\circ r(n)$.
Define $u(n)=t(n+1)$, $n\ge 1$.
We plan to show that $k_n=p(X)^{n+1}_n\circ h(n+1)_{n+1}$, $n\ge 1$,
induce a level morphism $k$ from $Z_u$ to $X$
such that $f\circ (k\circ p(Z)_u)=g$.
\par We need to prove $p(X)^{n+1}_n\circ k_{n+1}=k_n\circ 
p(Z)^{t(n+2)}_{t(n+1)}$
for all $n$.
If we show $f_{n+1}\circ k_{n+1}=f_{n+1}\circ h(n+1)_{n+1}\circ 
p(Z)^{t(n+2)}_{t(n+1)}$,
then it implies
$p(X)^{n+1}_n\circ k_{n+1}=p(X)^{n+1}_n\circ h(n+1)_{n+1}\circ 
p(Z)^{t(n+2)}_{t(n+1)}=
k_n\circ p(Z)^{t(n+2)}_{t(n+1)}$, i.e. what we need.
\par $f_{n+1}\circ h(n+1)_{n+1}\circ p(Z)^{t(n+2)}_{t(n+1)}=
h_{n+1}\circ r(n+1)_{n+1}\circ p(Z)^{t(n+2)}_{t(n+1)}=
h_{n+1}\circ p(Z)^{t(n+1)}_{s(n+1)}\circ p(Z)^{t(n+2)}_{t(n+1)}=
h_{n+1}\circ p(Z)^{t(n+2)}_{s(n+1)}$.
Also, $f_{n+1}\circ k_{n+1}=f_{n+1}\circ p(X)^{n+2}_{n+1}\circ h(n+2)_{n+2}=
p(Y)^{n+2}_{n+1}\circ f_{n+2}\circ h(n+2)_{n+2}=
p(Y)^{n+2}_{n+1}\circ h_{n+2}\circ r(n+2)_{n+2}=
p(Y)^{n+2}_{n+1}\circ h_{n+2}\circ p(Z)^{t(n+2)}_{s(n+2)}=
h_{n+1}\circ p(Z)^{s(n+2)}_{s(n+1)}\circ p(Z)^{t(n+2)}_{s(n+2)}=
h_{n+1}\circ p(Z)^{t(n+2)}_{s(n+1)}$. Thus
$f_{n+1}\circ k_{n+1}=f_{n+1}\circ h(n+1)_{n+1}\circ p(Z)^{t(n+2)}_{t(n+1)}$.
Also, we established
$f_{n}\circ k_{n}=
h_{n}\circ p(Z)^{t(n+1)}_{s(n)}$
which means $f\circ k=h\circ p(Z)^u_s$. Composing with $p(Z)_u$
gives $f\circ k\circ p(Z)_u=h\circ p(Z)_s=g$.
\par General Case. Using \ref{XXX1.8} and \ref{XXX1.mono} we can find
an increasing sequence $s:N\to N$ such that
$f'=p(X)_s\circ f:X_s\to Y$ is a level morphism
and for any two morphisms
$a,b:P\to X_{s(n+1)}$ of $C$ the equality $f'_{n+1}\circ a=f'_{n+1}\circ b$
implies $p(X_s)^{n+1}_n\circ a=p(X_s)^{n+1}_n\circ b$.
Notice that $p(X)_s:X\to X_s$ is an isomorphism
as $s(N)$ is cofinal in $N$.
By the Special Case
$f'_\ast:Mor(Z,X_s)\to Mor(Z,Y)$ is a bijection for all
sequentially movable objects $Z$ of pro-$C$.
As $p(X)_s$ is an isomorphism,
$f_\ast:Mor(Z,X)\to Mor(Z,Y)$ is a bijection for all
sequentially movable objects $Z$ of pro-$C$.
\end{pf}

\begin{Rem}\label{XXX3.15.5}
The above result is not valid for arbitrary monomorphisms
$f:X\to Y$. The pro-group $G$ in \ref{XXX3.11} is sequentially
movable, the trivial morphism $f:0\to G$
has the property that $f_\ast:Mor(P,0)\to Mor(P,G)$ is a bijection for all
groups $P$ but
$f_\ast:Mor(G,0)\to Mor(G,G)$ is not a surjection.
\end{Rem}

\medskip
\medskip

\section{Stability in pro-categories}

\medskip

The next two results will be useful for applications
to pro-groups.
\begin{Prop}\label{XXX3.16} If $X$ is a pro-object such that
each $p(X)^\beta_\alpha$ is a monomorphism of $C$,
then the following conditions are equivalent:

1. $X$ is stable.

2. $X$ is movable.

3. There is $\alpha\in I(X)$
such that each $p(X)^\beta_\alpha$,
$\beta>\alpha$, is an isomorphism of $C$.
\end{Prop}
\begin{pf} 1)$\implies$2) follows from the fact
that movability is preserved by isomorphisms
(see \cite{MS}, Theorem 1 on p.159).
\par 2)$\implies$3). Let $\alpha\in I(X)$ and pick $\beta>\alpha$
such that for each $\gamma>\beta$ there is
$r_\gamma:X_\beta\to X_\gamma$ with
$p(X)^\gamma_\alpha\circ r_\gamma=p(X)^\beta_\alpha$.
Now $p(X)^\beta_\alpha\circ (p(X)^\gamma_\beta\circ 
r_\gamma)=p(X)^\beta_\alpha\circ id(X_\beta)$, so
$p(X)^\gamma_\beta\circ r_\gamma=id(X_\beta)$.
That means $p(X)^\gamma_\beta$ has a left inverse and must be an isomorphism
by \ref{XXX1.2}.
\par 3)$\implies$1).
Notice that $p(X)_\alpha:X\to X_\alpha$
is an isomorphism.
\end{pf}

\begin{Prop}\label{XXX3.17} (Corollary 2.12 of \cite{DR2}) Suppose 
$X$ is a pro-object such that
each $p(X)^\beta_\alpha$ is an epimorphism of $C$.
$X$ is stable if and only if there is $\alpha\in I(X)$
such that $p(X)^\beta_\alpha$ is an isomorphism of $C$
for all $\beta>\alpha$.
\end{Prop}

In \cite{DR2} the authors discussed the question of
objects of $C$ having stable images (respectively, stable subobjects)
in pro-$C$. In this section we deal with analogous
question of objects having stable strong images (respectively, stable strong
subobjects). We will use repeatedly the following result.

\begin{Cor} (3.6 of \cite{DR2})
\label{XXX3.6ofDR2} Suppose $C$ is a balanced category with epimorphic images. 
If pro-$C$ is balanced, then for any epimorphism (respectively, monomorphism) 
$f:X\to Y$ of pro-$C$ there exists a level morphism $f':X'\to Y'$ and 
isomorphisms $i:X\to X'$, $j:Y'\to Y$ such that 
$f=j\circ f'\circ i$, $I(X')$ is a cofinite directed set, 
and $f'_\alpha$ is an epimorphism (respectively, monomorphism) of $C$ for each 
$\alpha\in I(Y')$. Moreover, the bonding morphisms of $X'$ 
(respectively, $Y'$) are chosen from the set of bonding morphisms of $X$ 
(respectively, $Y$). 
\end{Cor} 

Recall that $C$ is a balanced category with epimorphic images
if every morphism $f$ of $C$ has a unique, up to isomorphism,
decomposition $f=f''\circ f'$ such that $f'$ is an epimorphism of $C$
and $f''$ is a monomorphism of $C$.

\begin{Def}\label{XXX4.n0} Let $C$ be a category.
$Y\in Ob(pro-C)$ is called a {\bf strong image} (respectively,
{\bf a strong subobject}) of an object $X$ of pro-$C$
  provided there is
a strong epimorphism (respectively, strong monomorphism)
$f:X\to Y$ (respectively, $f:Y\to X$) of pro-$C$.
\end{Def}


\begin{Def}\label{XXX4.n1} Let $C$ be a category.
An object $P$ of $C$
has {\bf stable strong images} (respectively, {\bf stable strong subobjects})
if any strong image (respectively, strong subobject) $X\in Ob(pro-C)$ 
of $P$ is stable.
\end{Def}


\begin{Thm}\label{XXX4.n2} Let $R$ be a principal ideal domain.
If $P$ is a finitely generated $R$-module, then it
has stable strong images (respectively, stable strong subobjects)
in the pro-category pro-$M_R$ of the category $M_R$ of $R$-modules.
\end{Thm}
\begin{pf} Suppose $f:X\to P$ is a strong monomorphism.
Since $M_R$ is a balanced category with epimorphic
images and pro-$M_R$ is balanced, \ref{XXX3.6ofDR2}
allows us to reduce the proof to the case where
$f$ is a level morphism, and each $f_\alpha $ is a monomorphism.
In particular, as $f_\beta=f_\alpha \circ p(X)^\beta _\alpha $,
each $p(X)^\beta _\alpha $ is a monomorphism.
Subsequently, we may simply identify all $X_\alpha $ with submodules
of $P$ so that all $p(X)^\beta _\alpha $ are inclusion-induced.
Now it suffices to show that there is
$\alpha \in I(X)$ such that $X_\beta =X_\alpha $ for all $\beta >\alpha $.
\par Special Case. $P$ is a torsion $R$-module.
In this case $P$ satisfies the descending chain condition
on submodules in view of Theorem 1.5 in \cite{Hun} (p.373), so
$X$ is stable.
\par Notice that one has the functor $Tor:M_R\to M_R$
such that $Tor(Q)$ is the torsion part of an $R$-module $Q$.
That functor can be extended to $Tor:pro-M_R\to pro-M_R$.
Therefore $Tor(X)$ is a strong subobject of $Tor(P)$
and must be stable by Special Case. Without loss of generality
we may assume $Tor(X_\alpha )=Tor(X_\beta )$ for all $\beta >\alpha $.
Also, as the rank of the free part of $X_\alpha $ is at most
the rank of the free part of $P$, we may assume that
the ranks of all free parts of modules $X_\alpha $ are equal to
a fixed natural number $m$.
Suppose $X$ is not stable. By \ref{XXX1.6b} there is a triple
$\gamma>\beta >\alpha $ of elements of $I(X)$ such that
for some morphism $r:P\to X_\beta $
one has $r|X_\gamma$ is the inclusion and $X_\gamma\ne X_\beta \ne X_\alpha $.
Notice that $r|T$ is the inclusion, where $T$ is the torsion submodule
of $X_\alpha $, so one can
put $Q=P/T$, $Y_\delta =X_\delta /T$ for all $\delta \in I(X)$,
and get $s:Q\to Y_\beta $
such that $s|Y_\gamma$ is the inclusion and $Y_\gamma\ne Y_\beta \ne Y_\alpha $
are all free $R$-modules of the same rank.
Let $M=\{x\in Y_\alpha \mid s(x)=x\}$.
It is a proper submodule of $Y_\alpha $ containing $Y_\gamma$ and 
contained in $Y_\beta $,
so it has the same rank as $Y_\alpha $ and cannot be a direct summand 
of $Y_\alpha $.
Therefore $Y_\alpha /M$ is not torsion-free and there is $u\in 
Y_\alpha \setminus M$
with $q\cdot u\in M$ for some $q\in R\setminus\{0\}$. Thus $q\cdot (s(u)-u)=0$
implying $s(u)-u=0$ and $u\in M$, a contradiction.
\par Suppose $f:P\to X$ is a strong epimorphism.
Since $M_R$ is a balanced category with epimorphic
images and pro-$M_R$ is balanced, \ref{XXX3.6ofDR2}
allows us to reduce the proof to the case where
$f$ is a level morphism and each $f_\alpha $ is an epimorphism.
Since $f_\alpha =p(X)^\beta _\alpha \circ f_\beta $,
each $p(X)^\beta _\alpha $ is an epimorphism.
Now it suffices to show that there is
$\alpha \in I(X)$ such that $\ker(f_\beta) =\ker(f_\alpha) $ for all 
$\beta >\alpha $.
\par Special Case. $P$ is a torsion $R$-module.
In this case $P$ satisfies the descending chain condition
on submodules in view of Theorem 1.5 in \cite{Hun} (p.373), so
$X$ is stable.
\par Notice that $Tor(X)$ is a strong image of $Tor(P)$
and must be stable by Case 1. Without loss of generality
we may assume $p(X)^\beta _\alpha |Tor(X_\beta )$
sends $Tor(X_\beta )$ isomorphically onto $Tor(X_\alpha )$
  for all $\beta >\alpha $.
Also, as the rank of the free part of $X_\alpha $ is at most
the rank of the free part of $P$, we may assume that
the ranks of all free parts of modules $X_\alpha $ are equal to
a fixed natural number $m$.
Therefore $\ker(p(X)^\beta _\alpha )$ must be a torsion
module for all $\beta >\alpha $. Since
$p(X)^\beta _\alpha |Tor(X_\beta )$
sends $Tor(X_\beta )$ isomorphically onto $Tor(X_\alpha )$
  for all $\beta >\alpha $, $p(X)^\beta _\alpha$ must be
an isomorphism.
\end{pf}

The same proof as in \ref{XXX4.n2} yields the following.
One cannot derive \ref{XXX4.n3} formally from \ref{XXX4.n2}
as the category of groups is larger than the category of
$Z$-modules (i.e., the category of Abelian groups).


\begin{Cor}\label{XXX4.n3}
If $P$ is a finitely generated Abelian group, then it
has stable strong images (respectively, stable strong subobjects)
in the category of pro-groups.
\end{Cor}

The remainder of this section is devoted to describing
another class of Abelian groups whose members
have stable strong images (respectively, stable strong subobjects)
in the category of pro-groups.


\begin{Def}\label{XXX4.n4}(\cite{Fu},p.29)
Let $S$ be a subset of an Abelian group $G$.
$S$ is called {\bf linearly independent}, or briefly,
{\bf independent}, if any relation
$$n_1\cdot a_1+\ldots+n_k\cdot a_k=0$$
implies $n_i\cdot a_i=0$ for all $i$.
\end{Def}


\begin{Def}\label{XXX4.n5}(\cite{Fu},p.31)
Let $G$ be an Abelian group $G$.
The {\bf rank} $r(G)$ of $G$ is the cardinality of
  a maximal independent subset of $G$ whose elements have orders
being prime or infinity.
\end{Def}


\begin{Thm}\label{XXX4.n6}
Let $G$ be a divisible Abelian group.
  $G$
has stable strong images (respectively, stable strong subobjects)
in the category of pro-groups if and only if its rank is finite.
\end{Thm}
\begin{pf} Suppose the rank of $G$ is infinite.
By Theorem 19.1 of \cite{Fu} (p.64) one can express
$G$ as the direct sum of an infinite sequence $G_i$, $i>0$,
of some of its non-trivial subgroups.
Let $F_n$ be the direct sum of $G_i$, $i\leq n$. Notice that
one has a tower $F$ bonded by projections such that the inclusion
$F\to G$ is a strong monomorphism and $F$ is not stable.
Similarly, let $H_n$ be the direct sum of $G_i$, $i\ge n$. Notice that
one has a tower $H$ bonded by inclusions such that there is
a strong epimorphism $G\to H$ and $H$ is not stable.
\par
Assume the rank $r$ of $G$ is finite.
\par Claim: Given a descending chain $G_i$ of divisible subgroups
of $G$ there is $k$ such that $G_{i+1}=G_i$ for $i\ge k$.
\par Proof of Claim. It suffices to show that one cannot
have a descending sequence $G_0=G\supset G_1\supset\ldots\supset G_{r+1}\ne 0$
such that each $G_{i+1}$ is a proper divisible subgroup of $G_i$
for $i=0,\ldots,r$. By Theorem 18.1 of \cite{Fu} (p.62)
each $G_{i+1}$ is a direct summand of $G_i$. Therefore, starting from
a maximal independent subset $S_{r+1}$ of $G_{r+1}$
(consisting of elements whose orders are prime or infinity)
one can increase it to a maximal independent subset $S_{r}$ of $G_{r}$
(consisting of elements whose orders are prime or infinity)
and so on. In the end we will get a maximal independent subset $S_{0}$ of $G$
(consisting of elements whose orders are prime or infinity)
whose cardinality is larger than $r$, a contradiction.
\par
Suppose $f:X\to G$ is a strong monomorphism.
Since $Gr$ is a balanced category with epimorphic
images and pro-$Gr$ is balanced, \ref{XXX3.6ofDR2}
allows us to reduce the proof to the case where
each $p(X)^\beta _\alpha $ is a monomorphism,
$f$ is a level morphism, and each $f_\alpha $ is a monomorphism.
Subsequently, we may simply identify all $X_\alpha $ with submodules
of $P$ so that all $p(X)^\beta _\alpha $ are inclusion-induced.
Now it suffices to show that there is
$\alpha \in I(X)$ such that $X_\beta =X_\alpha $ for all $\beta >\alpha $.
Suppose $X$ is not stable. We may reduce the general case to the one
of $X$ being a tower such that for each $n$ there is
a homomorphism $r_n:G\to X_n $
so that $r_n|X_{n+1}=id$ and $X_{n+1}$ is a proper subgroup
of $X_n$.
Let $G_n=r_n(G)$. It is a divisible subgroup of $G$
and $X_{n+1}\subset G_n\subset X_n$.
There is $k$ such that $G_{n+1}=G_n$ for $n\ge k$
implying $X_{n+1}=X_n$ for $n>k$, a contradiction.
\par Suppose $f:G\to X$ is a strong epimorphism.
Since $Gr$ is a balanced category with epimorphic
images and pro-$Gr$ is balanced, \ref{XXX3.6ofDR2}
allows us to reduce the proof to the case where
each $p(X)^\beta _\alpha $ is an epimorphism,
$f$ is a level morphism, and each $f_\alpha $ is an epimorphism.
Now it suffices to show that there is
$\alpha \in I(X)$ such that $\ker(f_\beta) =\ker(f_\alpha) $ for all 
$\beta >\alpha $.
Suppose $X$ is not stable. Again, we can reduce the general case to that
of $X$ being a tower, $\ker(f_{n+1})$ being a proper
subgroup of $\ker(f_n)$ for each $n$, and the equality
$f_n\circ r_n=p(X)^{n+1}_n$ for some homomorphism
$r_n:X_{n+1}\to G$.
Moreover, we may assume that every divisible subgroup of $G$
contained in all $\ker(f_n)$, $n\ge 1$, is trivial. Indeed,
one can consider a maximal divisible subgroup $D$ of the intersection
of all $\ker(f_n)$, $n\ge 1$. $D$ is a direct summand of $G$,
so replacing $G$ by $G/D$ does the trick.
Let $H_n$ be the subgroup of $G$ generated by elements $x-r_mf_{m+1}(x)$,
where $x\in G$ and $m\ge n$. Notice that $H_n$ is divisible and 
$H_{n+1}\subset H_n$.
Let us show $H_n\subset \ker(f_n)$.
Indeed, $f_n(r_mf_{m+1})=(p(X)^m_n\circ f_m)\circ r_m\circ f_{m+1}=
p(X)^m_n\circ p(X)^{m+1}_m\circ f_{m+1}=f_n$ for $m\ge n$,
so $f_n(x-r_mf_{m+1}(x))=0$.
Pick $k\ge 1$ such that $H_m=H_n$ for $m,n\ge k$. That means $H_n=0$
for $n\ge k$ proving that $f:G\to X$ is an isomorphism.
Thus $X$ is stable, a contradiction.
\end{pf}

\medskip
\medskip

\section{Bimorphisms in pro-categories}

\medskip

The purpose of this section is to relate bimorphisms
of pro-$C$ to bimorphisms of tow($C$) for categories $C$
with direct sums and weak push-outs.
\begin{Lem}\label{XXX4.1} Suppose $C$ is a category with direct sums
and weak push-outs. If $f:X\to Y$ is a level morphism
of pro-$C$ which is a bimorphism of pro-$C$,
then the set of sequences $s$ in $I(X)$
such that $f_s:X_s\to Y_s$ is a bimorphism
of tow($C$) is cofinal among all sequences in $I(X)$.
\end{Lem}
\begin{pf} Let $g:I(X)\times I(X)\to I(X)$
be a function such that $g(\alpha,\beta)>\alpha,\beta$.
Using \ref{XXX1.mono} and \ref{XXX1.epi} there are functions $m:I(X)\to I(X)$
and $e:I(X)\to I(X)$ with the following properties:

1. $m(\alpha)>\alpha$ and for any two
morphisms $a,b:P\to X_{m(\alpha)}$ the equality
$f_{m(\alpha)}\circ a=f_{m(\alpha)}\circ b$ implies
$p(X)^{m(\alpha)}_\alpha\circ a=p(X)^{m(\alpha)}_\alpha\circ b$.

2. $e(\alpha)>\alpha$ and for any two
morphisms $a,b:Y_{\alpha}\to P$ the equality
$a\circ f_{\alpha}=b\circ f_{\alpha}$ implies
$a\circ p(Y)^{e(\alpha)}_\alpha=b\circ p(Y)^{e(\alpha)}_\alpha$.
\par Given any sequence $t$ in $I(X)$
define $s(1)=g(m(t(1)),e(t(1)))$ and, inductively,
$s(n+1)=g(g(m(t(n)),e(t(n))),s(n))$.
Using \ref{XXX1.mono} and \ref{XXX1.epi} it is easy to check that 
$f_s$ is a bimorphism
of pro-$C$.
\end{pf}

\begin{Thm}\label{XXX4.2} If $C$ is a category with direct sums
and weak push-outs, then the following conditions are equivalent:

1. tow($C$) is balanced.

2. pro-$C$ is balanced.
\end{Thm}
\begin{pf} 1)$\implies$2). Suppose $f:X\to Y$
is a level morphism of pro-$C$ which is a bimorphism of pro-$C$.
We can find a cofinal subset $\Sigma$ of the set of
increasing sequences in $I(X)$ such that $f_s:X_s\to Y_s$
is a bimorphism of pro-$C$ for each $s\in \Sigma$.
Now, each $f_s$ is an isomorphism, so $f$ is an isomorphism
by \ref{XXX1.11}.
\par 2)$\implies$1). This amounts to showing that any
bimorphism of tow($C$) is also a bimorphism of pro-$C$. That was done 
in \ref{XXX1.10}.
\end{pf}

\medskip
\medskip

\section{Weak equivalences in pro-homotopy}

\medskip

Recall that {\bf a weak equivalence} in pro-$H_0$
is a morphism $f:X\to Y$ such that pro-$\pi_n(f)$
is an isomorphism for all $n$.
Also, {\bf the deformation dimension} $\dim_{def}(X)$
of $X$ is the smallest number $n$ such that for any
$\alpha\in I(X)$ there is $\beta>\alpha$
with $p(X)^\beta_\alpha$ having a representative
with image contained in the $n$-skeleton of $X_\alpha$
(see \cite{D1}).
\par The purpose of this section is to generalize
some versions of the Whitehead Theorem in pro-homotopy
(see \cite{D1} and \cite{MS}).

\begin{Lem}\label{XXX5.0} Any weak equivalence $g:X\to Y$ of tow($H_0)$
has the property that $g_\ast:Mor(P,X)\to Mor(P,Y)$
is surjective for all CW complexes $P$.
\end{Lem}
\begin{pf} \par Notice that every object
of tow$(H_0)$ is equivalent to
a tower
of spaces homotopically equivalent to pointed connected CW complexes
so that bonding maps are Hurewicz fibrations.
(see \cite{D1}, Theorem 5.2 and its proof).
Let us assume $X$ and $Y$ are towers in the category
of spaces homotopically equivalent to pointed connected CW complexes
as objects and Hurewicz fibrations as morphisms.
Using \ref{XXX1.8} we can reduce the proof to the case of $f$
being a level morphism. Moreover, as $p(Y)^m_n$
are Hurewicz fibrations, we may assume that $f_n$
are actually maps (as opposed to homotopy classes of maps)
so that $p(Y)^m_n\circ f_m=f_n\circ p(X)^m_n$ for $m>n$.
Let $\bar X$ (respectively, $\bar Y$) be the inverse limit
of $X$ (respectively, $Y$) and let $\bar f:\bar X\to \bar Y$
be the map induced by $f$.
Notice that $Mor(P,\bar X)\to Mor(P,X)$ is an epimorphism
for all pointed connected CW complexes $P$, and the same statement
holds for $Y$. It follows from the fact that bonding maps are
Hurewicz fibrations (see \cite{D1}, Theorem 5.2 and its proof).
Therefore, it suffices to show that $\bar f$ is a weak homotopy
equivalence. By \cite{Bou-K} (p.254) one has the following
short exact sequence:
\begin{center}
\(
0\to \underset{\leftarrow}{\text{lim}}^1\pi_{i+1}(X)\to \pi_i(\bar 
X)\to \underset{\leftarrow}{\text{lim}}\ \pi_{i}(X)\to 0.
\)
\end{center}
Since the same sequence holds for $Y$, the Five Lemma implies that
$\bar f$ is a weak homotopy
equivalence.
\end{pf}

\begin{Thm}\label{XXX5.1} Suppose $f:X \to Y$ is a weak equivalence of pro-$H_0$.
If $Y$ is sequentially movable, then $f$ is
a strong epimorphism of pro-$H_0$.
\end{Thm}
\begin{pf} Assume that $f$ is a level morphism
of pro-$H_0$. As in \cite{D1} or using the same technique
as in the proof of \ref{XXX4.1} (see also \cite{MS}, p.160) one notices
that the set of sequences $s$ in $I(X)$ such that
$f_s$ is a weak equivalence is cofinal among
all sequences in $I(X)$.
By \ref{XXX3.14} and \ref{XXX5.0}, $f$ is a strong epimorphism of pro-$H_0$.
\end{pf}

\begin{Thm}\label{XXX5.2} Suppose $f:X \to Y$ is a weak equivalence of pro-$H_0$.
  $f$ is
an isomorphism of pro-$H_0$ in the following two cases:

1. $\dim_{def}(X)$ is finite and $Y$ is sequentially movable.

2. $\dim_{def}(Y)$ is finite and $X$ is sequentially movable.
\end{Thm}
\begin{pf} In case of 1) $f$ is a strong epimorphism
of pro-$H_0$ and it is shown in \cite{D1}
that $\dim_{def}(Y)\leq \dim_{def}(X)$ in such a case.
Therefore both $X$ and $Y$ are of finite deformation
dimension and $f$ is an isomorphism of pro-$H_0$ by
\cite{D1} (see also \cite{MS}, Theorem 3 on p.149).
\par In case of 2) there is a right inverse $g:Y\to X$
as shown in \cite{D1} (see \cite{MS}, Theorem 4 on pp.149--150). By 
case 1) $g$ is an isomorphism,
so $f$ is an isomorphism of pro-$H_0$ as well.
\end{pf}

\medskip
\medskip

\section{Bimorphisms in pro-homotopy}

\medskip

In this section we give partial answers to the following
question.
\begin{Problem}\label{XXX6.1} If \(f:X \rightarrow Y\) is a 
bimorphism in pro-\(H_0\), is \(f\) an isomorphism?
\end{Problem}

Notice that \(H_0\) has direct sums in the form of the wedge
of CW complexes.
Also, \(H_0\) has weak push-outs in the form of the union of mapping cylinders.

\begin{Prop}\label{XXX6.2} Suppose $f:X\to Y$ is a level morphism of pro-$H_0$
such that for every $\alpha\in I(X)$ there is $\beta>\alpha$
with the property that for any morphisms $a,b:\Sigma(P)\to X_\beta$
of $H_0$, the equality $f_\beta\circ a=f_\beta\circ b$
implies $p(X)^\beta_\alpha \circ a=p(X)^\beta_\alpha\circ b$.
If $f$ is an epimorphism of pro-$H_0$,
then it is a weak equivalence of pro-$H_0$.
\end{Prop}
  \begin{pf}
In case of $f$ being a morphism between towers
it follows from Theorem 2.10 of \cite{DR}.
Indeed, the above condition implies that $\pi_k(f)$
is a monomorphism for all $k\ge 1$
and 2.10 of \cite{DR} says that any epimorphism of tow($H_0)$
is a weak equivalence in such a case.
In the general case one reduces the problem
to $f_s:X_s\to Y_s$, where $s$ is a sequence in $I(X)$
so that $f_s$ satisfies the assumptions of this proposition.
The set of such $s$ is cofinal among all sequences in $I(X)$.
Since each $f_s$ is a weak equivalence, so is $f$.
\end{pf}

\begin{Thm}\label{XXX6.3}
If \(f: X \rightarrow Y\) is a bimorphism in pro-\(H_0\), then it is 
a weak equivalence.
\end{Thm}
\begin{pf} Assume $f$ is a level morphism. Use \ref{XXX6.2} and 
\ref{XXX1.mono}.
\end{pf}

\begin{Thm}\label{XXX6.4}
Suppose  \(f: X \rightarrow Y\) is a bimorphism in pro-\(H_0\). Then 
\(f\) is an isomorphism
if one of the following conditions is satisfied:

i. \(Y\) is sequentially movable,

ii. $\dim_{def}(Y)$ is finite.
\end{Thm}
\begin{pf} By \ref{XXX6.3}, $f$ is a weak equivalence. In case of i) 
$f$ is a strong epimorphism by \ref{XXX5.1}.
Therefore, by \ref{XXX2.12}, it is an isomorphism.
\par In case of ii) $f$ has a right inverse
(see \cite{D1}), so it must be an isomorphism.
\end{pf}

\begin{Thm}\label{XXX6.5}
If \(f:X \rightarrow Y\) is a bimorphism in pro-\(H_0\), then 
$f_\ast:Mor(Z,X)\to Mor(Z,Y)$ is a bijection
for all sequentially movable objects $Z$  of pro-\(H_0\).
\end{Thm}
\begin{pf} Using \ref{XXX1.11} and as in \ref{XXX4.1} one can reduce 
it to $f$ being a level morphism
of tow($H_0$). Since $f$ is a weak equivalence,
then $f_\ast:Mor(P,X)\to Mor(P,Y)$ is a surjection
for all CW complexes $P$ (see \ref{XXX5.0}). By \ref{XXX3.15},
$f_\ast:Mor(Z,X)\to Mor(Z,Y)$ is a bijection
for all sequentially movable objects $Z$  of pro-\(H_0\).
\end{pf}

Recall that $\pi_\ast(P)$ is the group of homotopy classes
of maps from $\bigvee\limits_{i=1}^\infty S^i$ to $P$.
As a consequence of the above theorem one gets,
in view of \ref{XXX3.3}, the following.
\begin{Cor}\label{XXX6.6}
If \(f: X \rightarrow Y\) is a bimorphism in pro-\(H_0\) and 
\(pro-\pi_*(Y)\) is
uniformly movable
then \(f_*:pro-\pi_*(X) \rightarrow pro-\pi_*(Y)\) is an isomorphism.
\end{Cor}
\begin{pf}
Notice that $f_\ast$ is a monomorphism for any bimorphism
$f:X\to Y$ by part b) of \ref{XXX1.mono} applied to
$P=\bigvee\limits_{i=1}^\infty S^i$.
\par Let $g=\underset{\leftarrow}{\text{lim}}(f_\ast)$.
Applying  \ref{XXX6.5} to $Z=\bigvee\limits_{i=1}^\infty S^i$
one gets that $g$ is an isomorphism as 
$\underset{\leftarrow}{\text{lim}}(pro-\pi_\ast(X))$ equals
$Mor_{pro-H_0}(Z,X)$ and 
$\underset{\leftarrow}{\text{lim}}(pro-\pi_\ast(Y))=Mor_{pro-H_0}(Z,Y)$.
Now, Part c) of \ref{XXX3.3} says that $f_\ast$ is an isomorphism
if pro-$\pi_\ast(Y)$ is uniformly movable.
\end{pf}
\medskip
\medskip

\section{Bimorphisms in the shape category}

\medskip

Let $HT_0$ be the homotopy category of pointed connected topological
spaces. A {\bf shape system} of $X\in Ob(HT_0)$
is an object $K$ of pro-$H_0$
such that for some morphism
$f:X\to K$ of pro-$HT_0$ the induced function
$f^\ast:Mor(K,L)\to Mor(X,L)$ is a bijection for
all $L\in Ob(H_0)$.
There is (see \cite{MS}) a {\bf shape category} $Sh$ and the shape functor
$S:Sh\to pro-H_0$ such that $S(X)$ is the shape system
for each pointed connected topological space $X$
and $S$ establishes a bijection between $Mor_{Sh}(X,Y)$
and $Mor_{pro-H_0}(S(X),S(Y))$.
In this sense one can identify $X$ with $S(X)$
and consider $Sh$ to be the full subcategory of pro-$H_0$
whose objects are shape systems of pointed connected topological
spaces. This is the approach we take in this section.

\begin{Def}\label{XXX7.1}
Given a directed set \(A\) and pointed CW-complexes 
\(\{P_{\alpha}\}_{\alpha \in A}\)
let \(WC(\{P_{\alpha}\}_{\alpha \in A})\) be the topological
space with the underlying set
  \(\bigvee\limits_{\alpha \in A} Cone(P_{\alpha})\) (we denote the 
base point of it by \(p\))
so that a set \(U\) is open if and only if the following conditions are satisfied:

1. \(U\cap Cone(P_{\alpha})\) is open in \(Cone(P_{\alpha})\) for 
each \(\alpha\).

2. If \(p\in U\), then there is \(\alpha_0 \in A\) such that
\(P_{\beta} \subset U\) for all \(\beta \geq \alpha_0\).
\end{Def}

\begin{Prop}\label{XXX7.2}
The space \(X=WC(\{P_{\alpha}\}_{\alpha \in A})\) is paracompact. 
A shape system for
\(WC(\{P_{\alpha}\}_{\alpha \in A})\) is
\((\bigvee\limits_{\beta \geq \alpha} \Sigma(P_{\beta}), 
j_{\alpha}^{\alpha'}, A)\), where \(j_{\alpha}^{\alpha'}\) is the 
natural inclusion.
\end{Prop}
\begin{pf} For each \(\alpha \in A\) let \(K_{\alpha}\) be the wedge
of all cones \(Cone(P_{\beta})\), where \(\beta\) is not bigger than 
or equal \(\alpha\), and
all suspensions \(\Sigma(P_{\beta})\) with \(\beta\ge\alpha\).
Let $\pi_{\alpha}:X\to K_{\alpha}$ be the projection so that $P_{\beta}$
(the base of the $Cone(P_{\beta})$)
is mapped to the base point for \(\beta\ge\alpha\). It is continuous
by the following argument: Suppose $V$ is an open subset of $K_{\alpha}$
and put
$U=\pi_{\alpha}^{-1}(V)$. $U\cap Cone(P_{\gamma})$ is open 
in $Cone(P_{\gamma})$ for all $\gamma$ as the projection
$Cone(P_{\gamma})\to \Sigma(P_{\gamma})$ is continuous for all $\gamma$.
If $p\in U$, then $V$ contains the base point of $K_{\alpha}$
and $U$ contains $P_{\beta}$ for all $\beta\ge \alpha$.
\par To show that $X$ is paracompact, it is sufficient to prove that,
for any open cover $\{U_s\}_{s\in S}$ of $X$ there is $\alpha\in A$
and an open cover $\{V_s\}_{s\in S}$ of $K_{\alpha}$
such that $\pi_{\alpha}^{-1}(V_s)\subset U_s$ for each $s\in S$.
Pick $s(0)\in S$ with $p\in U_{s(0)}$ and choose $\alpha\in A$
with $P_\beta\subset U_{s(0)}$ for $\beta\ge\alpha$.
Define $W_{s(0)}=U_{s(0)}$ and
 $W_s=U_s\setminus(\{p\}\cup \bigcup\limits_{\beta \ge\alpha } P_\beta )$ 
for $s\ne s(0)$.
Notice that for each $s\in S$ there is an open subset $V_s$ of
$K_\alpha$ such that $\pi_{\alpha}^{-1}(V_s)=W_s$.
Since $\{W_s\}_{s\in S}$ is an open cover of $X$, $\{V_s\}_{s\in S}$
is an open cover of $K_{\alpha}$ and the proof of $X$ being paracompact
is completed.
\par To prove the second part of the proposition it suffices
to show that \((K_{\alpha}, i_{\alpha}^{\alpha'}, A)\), where 
\(i_{\alpha}^{\alpha'}\) is the natural projection, is the shape 
system of $X$.
Since $i_{\alpha}^{\alpha'}\circ \pi_{\alpha'}=\pi_{\alpha}$
for $\alpha\leq \alpha'$, it suffices to show that the following
two statements are valid:
\par
a. Given any map $f:X\to K$ from $X$ to a CW complex $K$,
there is $\alpha\in A$ and a map $g:K_\alpha\to K$
such that $g\circ \pi_{\alpha}$
is homotopic to $f$.
\par
b. Given $\alpha\in A$ and given two maps $f,g:K_\alpha\to K$ from
$K_\alpha$ to a CW complex $K$ such that $f\circ \pi_{\alpha}$
is homotopic to
$g\circ \pi_{\alpha}$,
there is $\beta\ge\alpha$
such that $f\circ i_{\alpha}^\beta$
is homotopic to $g\circ i_{\alpha}^\beta$.
\par
Since every map to a CW complex is homotopic
to a locally compact map (see \cite{D3}), we will reduce a) and b) to the case
of $f$ and $g$ being locally compact.
\par Suppose $f:X\to K$ is a locally compact map from $X$
to a CW complex $K$. Let $C$ be a closed neighborhood of $p$ in $X$
such that $f(C)$ is contained in a compact subcomplex $L$ of $K$
containing the base point $\ast$ of $K$.
Since $L$ is locally contractible, there is a closed neighborhood
$D$ of $p$ in $X$ such that $f|D$ is homotopic to the constant
map. By the Homotopy Extension Theorem for locally compact maps
(see \cite{D3}),
$f$ is homotopic to $h:X\to K$ such that
$h(D)=\ast$. As in the proof of paracompactness of $X$,
there is $\alpha\in A$ and $g:K_\alpha\to K$ such that
$h=g\circ\pi_\alpha$.
\par Suppose $\alpha\in A$ and suppose $f,g:K_\alpha\to K$
are two locally compact maps from
$K_\alpha$ to a CW complex $K$ such that $f\circ \pi_{\alpha}$
is homotopic to
$g\circ \pi_{\alpha}$. As above, we may assume that both $f$ and $g$
are constant on some neighborhood of the basepoint of $K_\alpha$.
Also, we may assume that the homotopy $H$ joining
$f\circ \pi_{\alpha}$ and $g\circ \pi_{\alpha}$
is locally compact. By adjusting $H$ we can make it constant
on a neighborhood $U$ of $p$. Find $\beta\ge\alpha$ such that
$P_\gamma\subset U$ for all $\gamma\ge\beta$.
As above, $H$ can be factored through $K_\beta\times I$
which gives a homotopy joining
$f\circ i_{\alpha}^\beta$
and $g\circ i_{\alpha}^\beta$.
\end{pf}

Let us show that representatives of bimorphisms
of the shape category have the property as in \ref{XXX6.2}.

\begin{Prop}\label{XXX7.3aa}
Let \(f:X \rightarrow Y\) be a bimorphism in the shape category of pointed
connected topological spaces. If $f$ is represented
by a level morphism $g:S(X)\to S(Y)$ of shape systems of $X$
and $Y$, then $g$ is an epimorphism of pro-$H_0$ 
and for every $\alpha\in I(S(X))$ there is $\beta>\alpha$
with the property that for any morphisms $a,b:\Sigma(P)\to S(X)_\beta$
of $H_0$, the equality $f_\beta\circ a=f_\beta\circ b$
implies $p(S(X))^\beta_\alpha\circ a=p(S(X))^\beta_\alpha\circ b$.
\end{Prop}
\begin{pf} Let $D$ be the full subcategory of pro-$H_0$ whose objects
are shape systems of pointed connected topological spaces.
It is clear that $g$ is a bimorphism of $D$.
Notice that $g$ is an epimorphism
of pro-$H_0$. Indeed,
if $u,v:Y\to Z$ satisfy $ug=vg$,
then $(p(Z)_\alpha u)g=(p(Z)_\alpha v)g$ for all $\alpha \in I(Z)$.
Now, $Z_\alpha $ is an object of $D$, so 
$p(Z)_\alpha u=p(Z)_\alpha v$ for all $\alpha \in I(Z)$
which is the same as $u=v$.
It remains to show that for every $\alpha\in I(S(X))$ there is $\beta>\alpha$
with the property that for any morphisms $a,b:\Sigma(P)\to S(X)_\beta$
of $H_0$, the equality $f_\beta\circ a=f_\beta\circ b$
implies $p(S(X))^\beta_\alpha\circ a=p(S(X))^\beta_\alpha\circ b$.
If that property does not hold, then,
as in the proof of \ref{XXX1.mono} in \cite{DR2},
there is a system \(Z=(\bigvee\limits_{\beta \geq \alpha} 
\Sigma(P_{\beta}), j_{\alpha}^{\alpha'}, A)\), where 
\(j_{\alpha}^{\alpha'}\) is the natural inclusion,
such that for some morphisms $u,v:Z\to S(X)$
one has $u\circ g=v\circ g$ but $u\ne v$.
Since, by \ref{XXX7.2}, $Z$ is an object of $D$, one arrives
at a contradiction.
\end{pf}

In view of \ref{XXX6.2} one gets the following:
\begin{Thm}\label{XXX7.3}
If \(f:X \rightarrow Y\) is a bimorphism in the shape category of pointed
connected topological spaces, then \(f\) is a weak equivalence.
\end{Thm}

\begin{Prop}\label{XXX7.4}
  Let $X$ be a pointed connected space.
If $(K_\alpha,p^\beta_\alpha,A)$ is a shape system
of $X$,
then $(\Sigma K_\alpha,\Sigma p^\beta_\alpha,A)$ is a shape system
of $\Sigma X$.
\end{Prop}
\begin{pf} Given $g:\Sigma X\to P\in ANR$
one has the adjoint map $g'$ from $X$ to the loop space $\Omega P$.
Also, given $g:X\to \Omega P$ one has the adjoint map
from $\Sigma X$ to $P$ which will be denoted by $g'$ as well.
\par
Given $g:\Sigma X\to P\in ANR$
the adjoint $g':X\to \Omega P$ factors
as $g'\sim h'\circ p_\alpha$ for some
$h':K_\alpha\to \Omega P$. Now
$g\sim h\circ \Sigma p_\alpha$,
where $h:\Sigma K_\alpha\to P$ equals $(h')'$.
\par If $g,h:\Sigma K_\alpha\to P\in ANR$
so that $g\circ\Sigma p_\alpha\sim h\circ\Sigma p_\alpha$,
then $g'\circ p_\alpha\sim h'\circ p_\alpha$
and there is $\beta>\alpha$
with $g'\circ p^\beta_\alpha\sim h'\circ p^\beta_\alpha$,
i.e. $g\circ \Sigma p^\beta_\alpha\sim h\circ \Sigma p^\beta_\alpha$.
\end{pf}

\begin{Thm}\label{XXX7.5}
If \(f:X \rightarrow Y\) is a bimorphism in the shape category of 
pointed topological spaces, then $f_\ast:Mor(Z,X)\to Mor(Z,Y)$ is a 
bijection
for all sequentially movable spaces $Z$ which are suspensions
of some space $Z'$.
\end{Thm}
\begin{pf} Almost the same as in \ref{XXX6.5} if one uses \ref{XXX7.3aa}.
\end{pf}

\begin{Thm}\label{XXX7.6}
Suppose \(f:X \rightarrow Y\) is a bimorphism in the shape category 
of pointed topological spaces. If \(Y\) is sequentially movable, then 
\(f\) is an isomorphism
in the following two cases:

1. $Y$ is the suspension of a space $Y'$.

2. $X$ is the suspension of a space $X'$.
\end{Thm}
\begin{pf} 1. \ref{XXX7.5} implies the existence of a left inverse of $f$,
so $f$ is an isomorphism by \ref{XXX1.2}.
\par 2.
$f$ is a weak equivalence by \ref{XXX7.3}, so \ref{XXX5.1} says that
it is a strong epimorphism. Assume $f$ is a level morphism
of pro-$H_0$. Given $\alpha\in I(X)$ we can find $\beta>\alpha$
such that for any $a,b:\Sigma(P)\to X_\beta$ the condition
$f_\beta\circ a=f_\beta\circ b$ implies $p(X)_\alpha^\beta\circ a
=p(X)_\alpha^\beta\circ b$ (see \ref{XXX7.3aa}).
Since $f$ is a strong epimorphism, there is $\gamma>\beta$
and $r:Y_\gamma\to X_\beta$ such that
$f_\beta\circ r=p(Y)_\beta^\gamma$.
Now, $f_\beta\circ (r\circ f_\gamma)=
p(Y)_\beta^\gamma\circ f_\gamma= f_\beta\circ p(X)_\beta^\gamma$,
so $p(X)_\alpha^\beta\circ r\circ f_\gamma=
p(X)_\alpha^\beta\circ p(X)_\beta^\gamma=p(X)_\alpha^\gamma$
which proves that $f$ is a strong monomorphism as well.
By \ref{XXX2.10}, $f$ is an isomorphism.
\end{pf}

\begin{Problem}\label{XXX7.7}
  If \(f:X \rightarrow Y\) is a bimorphism in the shape category of 
pointed topological spaces,
is \(f\) a bimorphism in pro-\(H_0\)?
\end{Problem}

\begin{Problem}\label{XXX7.8}
  If \(f:X \rightarrow Y\) is a bimorphism of the shape category of 
pointed metric continua, is \(f\) a weak isomorphism? Is \(f\) an 
isomorphism?
\end{Problem}

\begin{Problem}\label{XXX7.9}
  If \(f:X \rightarrow Y\) is a bimorphism of the shape category of 
pointed movable metric continua,  is \(f\) a weak isomorphism?
\end{Problem}

\medskip
\medskip

\section{Borsuk's problem and strong monomorphisms}

\medskip

The following question comes up naturally.
\begin{Problem}\label{XXX8.0}
Let $P$ be a finite connected pointed CW complex.
Does $P$ have stable strong images (respectively, stable strong
subobjects) in pro-$H_0$?
\end{Problem}

In this section we give partial answers to the part of \ref{XXX8.0}
dealing with strong subobjects and we point out that it is stronger
than the following problem posed by K.Borsuk \cite{Bor}.
\begin{Problem}[K.Borsuk]\label{XXX8.1}
  Suppose $X_n$ is a sequence of compact ANRs such that
for each $n$ there is a retraction $r_n:X_n\to X_{n+1}$.
Is there a number $m$ such that all retractions $r_n$
are homotopy equivalences for $n>m$?
\end{Problem}

See \cite{Husch}, \cite{Mos}, \cite{San},
and \cite{Singh} for partial solutions to the above problem.

\par
If $X$ is the inverse limit of the inverse sequence
$(X_n,i^m_n)$, where $i^m_n$ is the inclusion for $m>n$,
then the above problem is equivalent to stability of $X$.
Indeed, in one direction it is quite obvious (if
$i^m_n$ are homotopy equivalences for $n$ large enough)
and in the other direction it follows that
$\pi_k(i^m_n)$ are isomorphisms for large $n$ and $k>\dim_{def}(X_1)$,
so $i^m_n$ must be homotopy equivalences.
\par First let us show how one creates a strong monomorphism
in the situation described by Borsuk's problem.
\begin{Prop}\label{XXX8.2} Suppose $C$ is a category and $X$
is an object of pro-$C$. If each $p(X)^{\beta}_\alpha$ has a left inverse
then there is a strong monomorphism
from $X$ to $P\in Ob(C)$.
\end{Prop}
\begin{pf} Pick $\gamma\in I(X)$. Given $\beta>\gamma$
let $r:X_\gamma\to X_\beta$ be the left inverse
of $p(X)^{\beta}_\gamma$.
Thus, $r\circ p(X)^{\beta}_\gamma=id(X_\beta )$.
This can be interpreted, in view of \ref{XXX2.5},
as a proof that $p(Y)_\gamma:Y\to Y_\gamma$
is a strong monomorphism, where $I(Y)=\{\beta\in I(X)\mid
\beta>\gamma\}$, $Y_\beta=X_\beta$, and $p(Y)^\sigma_\tau=
p(X)^\sigma_\tau$ for all $\sigma,\tau\in I(Y)$.
In other words, $Y$ is a subsystem of $X$ with $I(Y)$
cofinal in $I(X)$. Therefore $X$ and $Y$ are isomorphic
and $X$ admits a strong monomorphism to $X_\gamma$.
\end{pf}

\begin{Problem}\label{XXX8.3} Suppose  \(X\) is a pointed metric continuum.
Is \(X\) uniformly movable if it admits a strong monomorphism to a 
compact polyhedron?
\end{Problem}

\medskip

The last problem is stronger that Borsuk's one. The result
below provides the justification of it. Indeed,
every uniformly movable object $X$ of pro-$H_0$
admits a strong
  epimorphism \(Q \rightarrow X\) by \ref{XXX3.2}.


\begin{Cor} \label{XXX8.3a} Let $X$ be an object of pro-$H_0$.
If there exist polyhedra \(P,Q\), a monomorphism \(X \rightarrow P\), 
and a strong
  epimorphism \(Q \rightarrow X\), then \(X\) is stable.
\end{Cor}
\begin{pf} By \ref{XXX3.5}, $X$ is dominated by an object of $H_0$
and \cite{D1} or \cite{MS} (Theorem 4 on p.224) say that $X$ is stable.
\end{pf}

The spaces of the type \(WC(\{P_{\alpha}\}_{\alpha \in A})\)
 as in \ref{XXX7.1}, are 
not uniformly movable and they admit strong monomorphisms to 
non-compact polyhedra. On the other hand, if \(X\) is a pointed 
movable metric continuum and
there is a monomorphism \(X \rightarrow P \in Ob(H_0)\) in the shape 
category of pointed movable metric continua, then \(X\) is stable.

\par It is well-known that Borsuk's problem has positive
answer if $X_1$ is simply connected.
Let us give a positive solution to \ref{XXX8.3} if $\pi_1(P)$ is finite.

\begin{Thm}\label{XXX8.13} Suppose $f:X\to P$ is a strong monomorphism
of pro-$H_0$ such that $P$ is a compact CW complex.
If $pro-\pi_1(X)$ is pro-finite  or $\pi_1(P)$ is finite,
then $X$ is stable.
\end{Thm}
\begin{pf}
\par Assume that $f$ is a level morphism
induced by $\{f_\alpha:X_\alpha\to P\}_{\alpha\in I(X)}$.
Notice that the deformation dimension of $X$ is bounded
by $\dim(P)$ (see \cite{D1}).
In view of results in \cite{D1} it suffices to prove
that $\pi_n(X)$
is stable for all $n$.

Case 1: $\pi_1(P)$ is finite.
In this case $\pi_n(f):\pi_n(X)\to \pi_n(P)$
is a strong monomorphism of pro-$Gr$ and
$\pi_n(P)$
is finitely generated and Abelian if $n\ge 2$.
By \ref{XXX4.n3}, $\pi_n(X)$
is stable for all $n$.
\par Case 2: $\pi_1(X)$ is pro-finite.
By results of \cite{DR2} the pro-group $\pi_1(X)$ is stable
(since it is pro-finite and admits a monomorphism to a group).
Again, by results of \cite{DR2},
we may assume that
all $\pi_1(X_\alpha)$ are finite and
all $\pi_1(p(X)^\beta_\alpha)$ are
isomorphisms. Since $\pi_1(f)$
is a strong monomorphism,
we may assume that each $\pi_1(f_\alpha)$
has a left inverse.
Therefore we may think of $\pi_1(X)$ as a retract
of $\pi_1(P)$ and the kernel of the retraction $r$
can be killed by attaching 2-cells
along $r(a)\cdot a^{-1}$ for every generator
$a$ of $\pi_1(P)$. This way one gets
a finite CW complex $Q$ containing $P$ such that
the composition $X\to P\to Q$ is a strong
monomorphism and $\pi_1(Q)$ is finite.
By Case 1, $X$ is stable.
  \end{pf}

\begin{Cor}\label{XXX8.23} Suppose $X$ is an object of pro-$H_0$.
If there is a strong monomorphism $f:X\to P$ and an epimorphism
$g:Q\to X$ such that $P$ is a finite CW complex and $Q$ is a  CW complex,
then $X$ is stable.
\end{Cor}

A first step to solve \ref{XXX8.3} would be the following:

\begin{Problem}\label{XXX8.40} Let $X$ be an object of pro-$H_0$.
Is \(X\) stable if there exist polyhedra \(P,Q\), a strong 
monomorphism \(X \rightarrow P\), and an epimorphism \(Q \rightarrow 
X\)?
\end{Problem}

The following problem is important because of possible applications
to dynamical systems.
\begin{Problem}\label{XXX8.41}
Suppose $P$ be a finite polyhedron and $f:P\to P$ is a morphism of $H_0$.
Let $X$ be the tower in $H_0$ such that $X_n=P$ and $p(X)^{n+1}_n=f$ 
for each $n$.
Is \(X \) stable if it is uniformly movable?
\end{Problem}

\section{Appendix}

Let us prove categorical characterizations of strong monomorphisms
and strong epimorphisms mentioned earlier.
Notice that the category $Sets$ of sets and functions
is a category with direct sums, direct products, push-outs,
and pull-backs. Therefore its dual category $Sets^\ast$
has the same properties. Indeed, existence of push-outs (respectively, direct sums)
in the dual
category is equivalent to existence of pull-backs (respectively, direct
products) in the original 
category. 
\par First, we plan to show that both pro-$Sets$
and pro-$Sets^\ast$ have the property that every monomorphism
(respectively, epimorphism) is a strong monomorphism
(respectively, strong epimorphism).

\begin{Lem}\label{XXappe.1}
If 
\begin{center} \setlength{\unitlength}{1mm} 

\begin{picture}(63,28)

\put(5.3,22){\(X_\beta \)}

\put(58.3,22){\(Y_\beta  \)}

\put(5.3,0){\(X_\alpha\)}

\put(58.3,0){\(Y_\alpha \)}

\put(32,25){\(f_\beta \)}

\put(32,3){\(f_\alpha \)}

\put(61,11){\(p(Y)^\beta _\alpha \)}

\put(-5,11){\(p(X)^\beta _\alpha \)}

\put(15,23){\vector(1,0){37}}

\put(15,1){\vector(1,0){37}}

\put(6.3,19){\vector(0,-1){13}}

\put(59.3,19){\vector(0,-1){13}}

\end{picture} 
\end{center} \noindent
\par\noindent is a commutative diagram of $Sets$, then the following conditions
are equivalent:
\par 1. There is $r:Y_\beta \to X_\alpha $ such that
$f_\alpha \circ r=p(Y)^\beta _\alpha $.
\par 2. If $u,v:Y_\alpha \to P$ and $u\circ f_\alpha =v\circ f_\alpha $,
then $u\circ p(Y)^\beta _\alpha =v\circ p(Y)^\beta _\alpha $.
\end{Lem}
\begin{pf} 1)$\implies$2) is obvious.
\par 2)$\implies$1). Let $P=Y_\alpha /im(f_\alpha )$,
 $u:Y_\alpha \to P$ is the projection, and $v:Y_\alpha \to P$
is the constant map to the point of $P$ obtained by collapsing $im(f_\alpha )$.
Since $u\circ f_\alpha =v\circ f_\alpha $, 
$u\circ p(Y)^\beta _\alpha =v\circ p(Y)^\beta _\alpha $
which means exactly that $im(p(Y)^\beta _\alpha )\subset im(f_\alpha )$
in which case $r$ exists.
\end{pf}

\begin{Cor}\label{XXappe.2}
Every epimorphism of pro-$Sets$
is a strong epimorphism of pro-$Sets$.
\end{Cor}
\begin{pf} Assume $f=\{f_\alpha \}_{\alpha \in A}$
is a level morphism of pro-$Sets$ and pick $\alpha \in A$.
If $f$ is an epimorphism, then \ref{XXX1.epi} 
leads to a commutative diagram as in \ref{XXappe.1}.
By  \ref{XXX2.5}, $f$ is a strong epimorphism.
\end{pf}

\begin{Lem}\label{XXappe.3}
If 
\begin{center} \setlength{\unitlength}{1mm} 

\begin{picture}(63,28)

\put(5.3,22){\(X_\beta \)}

\put(58.3,22){\(Y_\beta  \)}

\put(5.3,0){\(X_\alpha\)}

\put(58.3,0){\(Y_\alpha \)}

\put(32,25){\(f_\beta \)}

\put(32,3){\(f_\alpha \)}

\put(61,11){\(p(Y)^\beta _\alpha \)}

\put(-5,11){\(p(X)^\beta _\alpha \)}

\put(15,23){\vector(1,0){37}}

\put(15,1){\vector(1,0){37}}

\put(6.3,19){\vector(0,-1){13}}

\put(59.3,19){\vector(0,-1){13}}

\end{picture} 
\end{center} \noindent
\par\noindent is a commutative diagram of $Sets$, then the following conditions
are equivalent:
\par 1. There is $r:Y_\beta \to X_\alpha $ such that
$r\circ f_\beta=p(X)^\beta _\alpha $.
\par 2. If $u,v:P\to X_\beta $ and $f_\beta  \circ u=f_\beta  \circ v$,
then $p(X)^\beta _\alpha\circ u =p(X)^\beta _\alpha\circ v$.
\end{Lem}
\begin{pf} 1)$\implies$2) is obvious.
\par 2)$\implies$1). Let $P$ be the point-inverse of
a point in $Y_\beta $ under $f_\beta $.
Let $u:P\to X_\beta $ be the inclusion and let
$v:P\to X_\beta $ be any constant function
to a point in $P$. Obviously, $f_\beta  \circ u=f_\beta  \circ v$,
so $p(X)^\beta _\alpha\circ u =p(X)^\beta _\alpha\circ v$
which means that $p(X)^\beta _\alpha (P)$ contains at most one
point. One defines $r:Y_\beta \to X_\alpha $
arbitrarily on $Y_\beta \setminus f_\beta (X_\beta )$
and $r(y)=p(X)^\beta _\alpha (f_\beta ^{-1}(y))$
for $y\in im(f_\beta )$.
\end{pf}

\begin{Cor}\label{XXappe.4}
Every monomorphism of pro-$Sets$
is a strong monomorphism of pro-$Sets$.
\end{Cor}
\begin{pf} Assume $f=\{f_\alpha \}_{\alpha \in A}$
is a level morphism of pro-$Sets$ and pick $\alpha \in A$.
If $f$ is a monomorphism, then \ref{XXX1.mono} 
leads to a commutative diagram as in \ref{XXappe.3}.
By  \ref{XXX2.5}, $f$ is a strong monomorphism.
\end{pf}

\begin{Cor}\label{XXappe.5}
Every monomorphism (respectively, epimorphism) of pro-$Sets^\ast$
is a strong monomorphism (respectively, strong epimorphism)
of pro-$Sets^\ast$.
\end{Cor}
\begin{pf}
Given a morphism $u:K\to L$ of $Sets^\ast$, we will denote
by $\mu(u):L\to K$ the corresponding function from $L$ to $K$.
Assume $f=\{f_\alpha \}_{\alpha \in A}$
is a level morphism of pro-$Sets^\ast$ and pick $\alpha \in A$.
\par Case 1. $f$ is a monomorphism of pro-$Sets^\ast$.
By  \ref{XXX1.mono} there is $\beta >\alpha $
with the property that, for any $u,v:P\to X_\beta $,
$f_\beta \circ u=f_\beta \circ u$
implies $p(X)^\beta _\alpha \circ u=p(X)^\beta _\alpha \circ v$.
That is the same as saying that for any functions
$u',v':X_\beta \to P$ the equality
$u'\circ \mu(f_\beta )=v'\circ \mu(f_\beta )$
implies $u'\circ \mu(p(X)^\beta _\alpha  )=v'\circ \mu(p(X)^\beta _\alpha )$.
By \ref{XXappe.1} there is a function $r:X_\alpha \to Y_\beta $
such that $\mu(f_\beta )\circ r=\mu(p(X)^\beta _\alpha )$.
If $g:Y_\beta \to X_\alpha $ is the morphism of $Sets^\ast$
corresponding to $r$, then $r\circ f_\beta =p(X)^\beta _\alpha $
which proves that $f$ is a strong monomorphism (see \ref{XXX2.5}).
\par The case of epimorphisms can be proved similarly
using \ref{XXappe.3}.
\end{pf}

\begin{Cor}\label{XXapp.1}
Suppose $f:X\to Y$ is a morphism of pro-$C$.
If, for any covariant functor $F:C\to D$, the induced morphism
$F(f):F(X)\to F(Y)$ is a monomorphism (respectively, epimorphism),
then $f$ is a strong monomorphism (respectively, strong epimorphism)
of pro-$C$.
\end{Cor}
\begin{pf} Assume $f=\{f_\alpha \}_{\alpha \in A}$
is a level morphism of pro-$C$ and pick $\alpha \in A$.
\par Case 1. $F(f)$ is a monomorphism for any covariant
functor $F:C\to D$. Consider $D=Sets^\ast$ and define $F(Z)=Mor_C(Z,X_\alpha )$
regarded as a covariant functor from $C$ to $D$.
Since $F(f)$ is a strong monomorphism by \ref{XXappe.5},
\ref{XXX2.5} says there is $\beta > \alpha $
and a morphism $r:F(Y_\beta )\to F(X_\alpha )$
such that $r\circ F(f_\beta )= F(p(X)^\beta _\alpha )$.
 That implies $p(X)^\beta _\alpha$ belongs to the image
of $F(f_\beta )$ (considered as a function
from $Mor_C(Y_\beta ,X_\alpha )$ to $Mor_C(X_\beta ,X_\alpha )$)
 and there is $g:Y_\beta \to X_\alpha $ with $g\circ f_\beta =
p(X)^\beta _\alpha $. Thus, see  \ref{XXX2.5}, $f$ is a strong monomorphism.
\par Case 2. $F(f)$ is an epimorphism for any covariant
functor $F:C\to D$. Consider $D=Sets$ and define $F(Z)=Mor_{pro-C}(Y,Z)$
regarded as a covariant functor from $C$ to $D$.
Since $F(f)$ is a strong epimorphism by \ref{XXappe.5},
\ref{XXX2.5} says there is $\beta > \alpha $
and a function $r:F(Y_\beta )\to F(X_\alpha )$
such that $F(f_\alpha)\circ r= F(p(Y)^\beta _\alpha )$.
Let $u=r(p(Y)_\beta )$.
Now, $f_\alpha \circ u=p(Y)^\beta _\alpha \circ p(Y)_\beta =p(Y)_\alpha $. 
Thus, see  \ref{XXX2.5}, $f$ is a strong epimorphism.
\end{pf}

\medskip
\medskip

\medskip
\medskip

      Jerzy Dydak,
Math Dept, University of Tennessee, Knoxville, TN 37996-1300. USA.
E-mail addresses: dydak@@math.utk.edu

\medskip
\medskip

      Francisco R.Ruiz del Portal,
Departamento de Geometr\'{\i}a y Topolog\'{\i}a,
Facultad de CC.Matem\'aticas.
Universidad Complutense de Madrid.
Madrid, 28040. Spain.
E-mail addresses: R\(_{-}\)Portal@@mat.ucm.es

\end{document}